\documentclass[12pt]{article}
\usepackage{amssymb}
\usepackage{eufrak}
\usepackage{amsmath}

\begin{document}
\newtheorem{thm}{Theorem}[section]
\newtheorem{lem}[thm]{Lemma}
\newtheorem{cor}[thm]{Corollary}
\newtheorem{prop}[thm]{Proposition}
\newtheorem{defin}[thm]{Definition}
\newtheorem{remark}[thm]{Remark}
\newtheorem{remarks}[thm]{Remarks}
\newtheorem{example}[thm]{Example}
\def\QED{\hbox{\hskip 1pt \vrule width4pt height 6pt depth 1.5pt \hskip 1pt}}
\newenvironment{prf}[1]{\trivlist
\item[\hskip \labelsep{\bf #1.\hspace*{.3em}}]}{~\hspace{\fill}~$\square$\endtrivlist}
\newenvironment{proof}{
\begin{prf}{Proof}}{
\end{prf}}
 \def\square{\QED}

\def\CX{{\mathbb C}}
\def\QX{{\mathbb Q}}
\def\NX{{\mathbb N}}
\def\ZX{{\mathbb Z}}
\def\DX{{\mathbb D}}
\def\AX{{\mathbb A}}
\def\ord{\mbox{ord }}
\def\GL{{\rm GL}}
\def\SL{{\rm SL}}
\def\gl{{\rm gl}}
\def\sl{{\rm sl}}
\def\M{{\rm gl}}
\def\d{{
\partial}}
\def\frakk{{\mathfrak k}}
\def\frakK{{\mathfrak K}}
\def\frakE{{\mathfrak E}}
\def\frakF{{\mathfrak F}}
\def\frakU{{\mathfrak U}}
\def\ofrakF{{\overline{\mathfrak F}}}
\def\oofrakF{{\overline{\mathfrak F}}_0}
\def\ofrakG{{\overline{\mathfrak G}}}
\def\oofrakG{{\overline{\mathfrak G}}_0}
\def\frakG{{\mathfrak G}}
\def\td{{\tilde{D}}}
\def\ttd{{\tilde{\tilde{D}}}}
\def \ldf{{{\rm LDF}_0}}
\def\QED{\hbox{\hskip 1pt \vrule width4pt height 6pt depth 1.5pt \hskip 1pt}}
\def\calD{{\cal D}}
\def\calC{{\cal C}}
\def\calS{{\cal S}}
\def\calT{{\cal T}}
\def\calF{{\cal F}}
\def\calE{{\cal E}}
\def\calL{{\cal L}}
\def\calG{{\cal G}}
\def\calK{{\cal K}}
\def\calO{{\cal O}}
\def\calU{{\cal U}}
\def\calM{{\cal M}}
\def\calP{{\cal P}}
\def\calA{\mbox{\sc a}}
\def\calB{{\cal B}}
\def\calCF{{\calC_0(\calF)}}
\def\gal{{\rm DiffGal}}
\def\C{{\rm Const}}
\def\d{{
\partial}}
\def\dx{{
\partial_{x}}}
\def\dt{{
\partial_{t}}}
\def\Autd{{{\rm Aut}_\Delta}}
\def\Gal{{\PGal}}
\def\PGal{{\rm Gal_\Delta}}
\def\Ga{{\bf G_a}}
\def\Gm{{\bf G_m}}
\def\ld{{l\d}}
\def\vr{{\Vec{r}}}
\def\vta{{\Vec{\tau}}}
\def\vt{{\Vec{t}}}
\def\KPV{{K^{\rm PV}_A}}
\def\ord{{\rm ord}}
\title{Galois Theory of Parameterized Differential Equations and Linear Differential Algebraic
Groups\footnote{This paper is an expanded version of a talk  presented at the conference {\em
Singularit\'es des  \'equations diff\'erentielles, syst\`emes int\'egrables et groupes quantiques},
November 24-27, 2004, Strasbourg, France.  The second author would like to thank the organizers of
this conference for inviting him.}}
\author{Phyllis J. Cassidy\footnote{The City College of New York, Department of Mathematics,
New York, New York 10038, USA pcassidy1@nyc.rr.com} \ and  Michael F. Singer\footnote{North Carolina
State University, Department of Mathematics, Box 8205, Raleigh, North Carolina 27695-8205, USA, 
singer@math.ncsu.edu. The second author was  supported by NSF Grant CCR- 0096842. }}

\date{February 14, 2005}
\maketitle
\begin{abstract} We present a Galois theory of parameterized linear differential equations where the
Galois groups are linear differential algebraic groups, that is, 
groups of matrices whose entries are functions of the parameters and satisfy a set of differential
equations with respect to these parameters. We present the basic constructions and results, give
examples, discuss how isomonodromic families fit into this theory and show how results from the
theory of linear differential algebraic groups may be used to classify systems of second order
linear differential equations.
\end{abstract}
\newpage
\tableofcontents
\newpage
\section{Introduction} We will describe a Galois theory of differential equations of the form
\[ \frac{\partial Y}{\partial x} = A(x,t_1, \ldots ,t_n)Y \]
where $A(x,t_1, \ldots ,t_n)$ is an $m \times m$ matrix with entries
 that are functions of the principal variable $x$ and parameters $t_1, \ldots ,t_n$.
The Galois groups in this theory are linear differential algebraic groups, that is, groups
of $m\times m$ matrices $(f_{i,j}(t_1, \ldots ,t_n))$ whose entries satisfy a fixed set of
differential equations.  For example, in this theory, the equation
\[ \frac{\partial y}{\partial x} = \frac{t}{x} y \]
has Galois group
\[ G = \{(f(t)) \ | \ f\frac{d^2f}{dt^2} - (\frac{df}{dt})^2 = 0 \} \ . \]
In the process, we will give an introduction to the theory of linear differential algebraic groups and show how one can use properties of
the structure to deduce results concerning parameterized linear differential equations.  \\[0.1in]
Various differential Galois theories now exist that go beyond the eponymous theory of linear differential equations
pioneered by Picard and Vessiot at the end of the $19^{th}$century and made rigorous and expanded
by Kolchin in the middle of the $20^{th}$ century. These include theories developed by B. Malgrange, A. Pillay, H. Umemura 
and one presently being developed by
P. Landesman. In many ways the Galois theory presented here is a special case of the results of Pillay and Landesman yet we hope that the 
explicit nature of our presentation and the applications we give justify  our exposition.  We will give a  comparison with these  theories
in the Final Comments.\\[0.1in]
The rest of the paper is organized as follows.  In section~\ref{sec1} we review the Picard-Vessiot theory of integrable systems of
linear partial differential equations.  In Section~\ref{sec2} we introduce and give the basic definitions and results for the Galois theory of
parameterized linear differential equations ending with a statement of the Fundamental Theorem of this  Galois theory as well as a characterization of
parameterized equations that are solvable in terms of parameterized liouvillian functions.  In Section~\ref{multgps} we describe the basic results concerning linear differential algebraic groups
and give many examples.  In Section~\ref{isomonsec} we show that, in the regular singular case, isomonodromic families 
of linear differential equations are precisely the parameterized linear differential equations whose parameterized 
Galois theory reduces to the usual Picard-Vessiot theory. In Section~\ref{secondsec} we apply a classification of $2 \times 2$ linear
differential algebraic groups to show that any parameterized system of linear differential equations with regular singular points is 
equivalent to a system that is
generic (in a suitable sense) or isomonodromic or solvable in terms of parameterized liouvillian functions.  Section~\ref{inversesec}
gives two simple examples illustrating the subtleties of the inverse problem in our setting.  In Section~\ref{finalsec} we discuss the relationship between the theory presented here
and other differential Galois theories and give some directions for future research.  The Appendices contain proofs of the results of Section~\ref{sec2}.

\section{Review of Picard-Vessiot Theory}\label{sec1} In the usual Galois theory of polynomial equations, the
Galois group is the collection of transformations of the roots that preserve all algebraic relations
among these roots.  To be more formal, given a field $k$ and a polynomial $p(y)$ with coefficients
in $k$, one forms the {\em splitting field} $K$ of $p(y)$ by adjoining all the roots of $p(y)$ to
$k$. The {\em Galois group} is then the group of all automorphisms of $K$ that leave each element of
$k$ fixed. The structure of the Galois group is well known to reflect the algebraic properties 
of the roots of $p(y)$. In this section we will review the Galois theory of linear differential
equations. Proofs (and  other references) can be found in \cite{PuSi2003}.\\
[0.1in]
One can proceed in an analogous fashion with integrable systems of linear differential equations and
define a Galois group that is a collection of transformations of  solutions of a linear differential
system that preserve all the algebraic relations among the solutions {\em and their derivatives}.
Let $k$ be a differential field\footnote{All fields in this paper will be of characteristic zero},
that is, a field $k$ together with a set of commuting derivations $\Delta = \{\d_1, \ldots ,\d_m\}$.
To emphasize the role of $\Delta$, we shall refer to such a field as a {\em $\Delta$-field}. 
Examples of such fields are the field ${\CX}(x_1, \ldots , x_m)$  of rational functions in $m$
variables, the 
quotient field $\CX((x_1, \ldots , x_m))$ of the ring of formal power series in $m$ variables  and
the quotient field $\CX(\{x_1, \ldots ,x_m\})$ of the ring of convergent power series, all  with the
derivations $\Delta =\{\frac{\d \ \ }{\d x_1},\ldots ,\frac{\d \ \ }{\d x_m}\}$. If $k$ is a
$\Delta$-field
and $\Delta' \subset \Delta$, the  field $C^{\Delta'}_k = \{ c\in k \ | \ \d c = 0 \mbox{ for all
} \d \in \Delta'\}$ is called the subfield of $\Delta'$-constants of $k$. When $\Delta' = \Delta$ we
shall write $C_k$ for $C^{\Delta}_k$ and refer to this later field as the field of constants of $k$.
An {\em integrable system of linear differential equations} is a set of equations
\begin{eqnarray}
\d_1 Y&=& A_1Y \label{intsys}\\
\d_2 Y  &=&  A_2 Y\notag\\
   &\vdots&  \notag\\
  \d_m Y  &= & A_m Y \notag
  \end{eqnarray}
  where the $A_i \in \M_n(k)$, the set of $n\times n$ matrices with entries in $k$, such that
 \begin{eqnarray}
\d_iA_j - \d_jA_i & = & [A_i, A_j] \label{intcond}
\end{eqnarray}
for all $i, j$.  These latter equations are referred to as the {\em integrability conditions}. Note
that if $m=1$, these conditions are trivially satisfied. \\
[0.1in]
The role of a splitting field is assumed  by  the {\em Picard-Vessiot extension} associated with the
integrable system (\ref{intsys}). This is a $\Delta$-extension field $K = k(z_{1,1}, \ldots ,
z_{n,n})$ where \begin{enumerate}
\item the $z_{i,j}$ are entries of a matrix $Z \in \GL_n(K)$ satisfying $\d_iZ = A_iZ$ for $i = 1,
\ldots ,m$, and
\item $C_K = C_k = C$, {\em i.e.,} the $\Delta$-constants of $K$ coincide with the  $\Delta$-constants of $k$.
\end{enumerate} Note that condition 1.~defines uniquely the actions on $K$ of the derivations $\d_i$
and that the integrability conditions (\ref{intcond}) must be satisfied since  these derivations
commute.
 We
refer to the $Z$ above as a {\em fundamental solution matrix} and we shall denote $K$ by $k(Z)$. If
$k =
\CX(x_1,
\ldots, x_m)$ with the obvious derivations, one can easily show the existence of Picard-Vessiot
extensions. If we let ${\Vec a}=(a_1, \ldots ,a_n)$ be a point of $\CX^n$ where the denominators of
all entries of the $A_i$ are holomorphic, then the Frobenius Theorem
(\cite{varadarajan_lie},Ch.~1.3) implies that, in a neighborhood of  ${\Vec a}$, there exist $n$
linearly independent analytic solutions $(z_{1,1}, \ldots , z_{n,1})^T, \ldots , (z_{1,n}, \ldots ,
z_{n,n})^T$ of the equations (\ref{intsys}). The field $k(z_{1,1}, \ldots ,z_{n,n})$ with the
obvious derivations satisfies the conditions defining a Picard-Vessiot extension. In general, if $k$
is an arbitrary $\Delta$-field with $C_k$ algebraically closed, then there always exists a
Picard-Vessiot extension $K$ for the integrable system (\ref{intsys}) and $K$ is unique up to
$k$-differential isomorphism. We shall refer to $K$ as the {\em PV-extension associated with
(\ref{intsys})}.
\\
[0.1in]
Let $K$ be a PV-extension associated with (\ref{intsys}) and let $K = k(Z)$ with $Z$ a fundamental
solution matrix.   If $U$ is another fundamental solution matrix then an easy calculation shows that
$\d_i(ZU^{-1}) = 0$ for all $i$ and so $Z U^{-1} \in \GL_n(C_k)$.  We define the {\em $\Delta$-Galois
group} $\Gal(K/k)$ of $K$ over $k$ (or of the system (\ref{intsys})) to be 
\[
\Gal(K/k) = \{\sigma :K \rightarrow K \ | \ \sigma \mbox{ is a $k$-automorphism of $K$ and }
\d_i \sigma = \sigma\d_i , \mbox{ for } i = 1 ,\ldots ,m\} \ .
\]
Note that a $k$-automorphism $\sigma$ of $K$ such that $\d_i\sigma = \sigma\d_i$ is called a $k$-{\em differential automorphism}.
For any $\sigma \in \Gal(K/k)$, we have that $\sigma(Z)$ is again a fundamental solution matrix so
the above discussion implies that $\sigma(Z) = ZA_{\sigma}$ for some $A_\sigma \in \GL_n(C_k)$. 
This yields a representation $\Gal(K/k) \rightarrow \GL_n(C_k)$.  Note that different fundamental
solution matrices yield conjugate representations.  A fundamental fact is that the image of
$\Gal(K/k)$ in $\GL_n(C_k)$ is Zariski-closed, that is, it is defined by a set of polynomial
equations involving the entries of the matrices and so has the structure of an linear algebraic
group.\\
[0.1in] These facts lead to a rich Galois theory, originally due to E.~Picard and E.~Vessiot and 
given rigor and greatly expanded by E.R.~Kolchin.  We summarize the fundamental result in the
following
\begin{thm}\label{PVthm} Let $k$ be a $\Delta$-field with algebraically closed field of constants
$C$ and (\ref{intsys}) be an integrable system of linear differential equations over $k$. 
\begin{enumerate}
\item There exists a PV-extension $K$ of $k$ associated with (\ref{intsys}) and this extension  is unique up to
differential $k$-isomorphism.
\item The $\Delta$-Galois group $\Gal(K/k)$ equals $G(C)$, where $G$ is a linear algebraic group defined over $C$.
\item The  map that sends any differential subfield $F, k \subset F \subset K$, to the group
$\Gal(K/F)$ is a bijection between the set of differential subfields of $K$ containing $k$ and and the set of algebraic subgroups  of $\Gal(K/k)$. Its inverse is given by the map that sends a Zariski closed group $H$ to
the field $K^H = \{z \in K \ | \ \sigma(z) = z \mbox{ for all } \sigma \in H\}$.
\item A Zariski-closed subgroup $H$ of $\PGal(K/k)$ is a normal subgroup of $\PGal(K/k)$ if and
only if the field $K^H$ is left set-wise invariant by   $\PGal(K/k)$. If this is the case, the 
map $\PGal(K/k) \rightarrow  \PGal(K^H/k)$ is surjective with kernel $H$ and
$K^H$ is a PV-extension of $k$ with PV-group isomorphic to $\PGal(K/k)/H$. Conversely, if $F$ is
a differential subfield of $K$ containing $k$ and $F$ is a PV-extension of $k$, then $\PGal(F/K)$
is a normal subgroup of $\PGal(K/k)$. 
\end{enumerate}
\end{thm}
\begin{remarks}{\em 1. The assumption that $C$ is algebraically closed is necessary for the
existence  of PV-extensions ({\em cf.,}\cite{seidenberg56}) as well as necessary to
guarantee that there are enough automorphisms so that 3.~is correct.  Kolchin's original development
of the Galois correspondence for PV-extensions did not make this assumption and he replaced
automorphisms of the PV-extension with embeddings of the PV-extension into a large field (a {\em
universal differential  field}) (see \cite{DAAG}). 
One can also study linear differential equations over fields whose fields of constants  are not
algebraically closed using descent techniques (see \cite{hoeij_put}). \\
[0.1in] 2. Theorem~\ref{PVthm} is usually stated and proven for the case when $m = 0$, the ordinary
differential case, although it is proven in this generality in \cite{DAAG}.  The usual proofs in
the ordinary differential case do
however usually generalize to this case as well. In the appendix of \cite{PuSi2003}, the authors
also discuss the case of $m>0$ and show how the Galois theory may be developed in this case. We will give a proof of a more general theorem in the appendix from which Theorem 2.1 follows as well.\\
[0.1in] 3. Theorem~\ref{PVthm} is a manifestation of a deeper result. If $K = k(Z)$ is a
PV-extension then the ring $k[Z, \frac{1}{\det Z}]$ is the coordinate ring of a torsor (principal
homogeneous space) $V$ defined over $k$ for the group $\Gal(K/k)$, that is, there is a morphism 
$V\times G
\rightarrow V$ denoted by $(v,g) \mapsto vg$ and defined over $k$ such that $v1 = v$ and $(vg_1)g_2
= v(g_1g_2)$ and such that the morphism $V\times G \rightarrow V\times V$ given by $(v,g) \mapsto
(v,vg)$ is an isomorphism. The path to the Galois theory given  by first establishing this fact is
presented in \cite{magid} and \cite{PuSi2003} (although Kolchin was well aware of this fact as
well,{\em cf.,} \cite{DAAG}, Ch.~VI.8 and the references there to the original papers.) This approach
allows one to
give an intrinsic definition of the linear algebraic group structure on the Galois group as well. 
}
\end{remarks} We end this section with a simple example that will also illuminate the Galois theory
of parameterized equations.
\begin{example} {\em Let $k = \CX(x)$ be the ordinary differential field with derivation
$\frac{d}{dx}$ and consider the differential equation
\[
\frac{dy}{dx} = \frac{t}{x}y
\]
where $t \in \CX$.  The associated Picard-Vessiot extension is $k(x^t)$. The Galois group will be
identified with a Zariski-closed subgroup of $\GL_1(\CX)$. When $t\in \QX$, one has that $x^t$ is an
algebraic function and when $t \notin \QX$, $x^t$ is transcendental.  It is therefore not surprising
that one can show that 
\begin{equation*}
\Gal(K/k) = 
\begin{cases}
\CX^* = \GL_1(\CX) & \text{if $t\notin \QX$,}\\
\ZX/q\ZX & \text{if $t = \frac{p}{q}, (p,q) = 1$.}
\end{cases}
\end{equation*} \hfill \QED}
\end{example}

\section{Parameterized Picard-Vessiot Theory}\label{sec2} In this section we will consider differential
equations of the form 
\[
\frac{\d Y}{\d x} = A(x,t_1, \ldots , t_m) Y
\]
where $A$ is an $n\times n$ matrix whose entries are functions of $x$ and parameters $t_1, \ldots ,
t_m$ and we will define a Galois group of transformations that preserves the algebraic relations
among a set of solutions and their derivatives {\em with respect to all the variables}. Before we
make things precise, let us consider an example.
\begin{example}{\em  Let $k = \CX(x,t)$ be the differential field with derivations $\Delta =\{\dx =
\frac{\d}{\d x}, \dt = \frac{\d}{\d t}\}$. Consider the differential equation
\[
\dx y = \frac{t}{x}y \ . 
\]
In the usual Picard-Vessiot theory, one forms the differential field generated by the entries of a
fundamental solution matrix and all their derivatives (in fact, because the matrix satisfies the
differential equation, we get the derivatives for free). We will proceed in a similar fashion here.
The function 
\[
y = x^t
\]
is a solution of the above equation. Although all derivatives with respect to $x$ lie in the field
$k(x^t)$, this is not true for $\dt(x^t) = (\log x) x^t$. Nonetheless, this is all that is missing
and  the derivations $\Delta$ naturally extend to the field 
\[
K = k(x^t, \log x),
\]
the field gotten by adjoining to $k$ a fundamental solution and its derivatives (of all orders) with
respect to all the variables. \\
[0.1in]
Let us now calculate the group $\Gal(K/k)$ of  $k$-automorphisms of $K$ commuting with both $\dx$ and
$\dt$. Let $\sigma \in \Gal(K/k)$.  We first note that $\d_x(\sigma(x^t) (x^t)^{-1}) = 0$ so $\sigma(x^t)
=
a_{\sigma} x^t$ for some $a_{\sigma} \in K$ with $\dx a_{\sigma} = 0$, {\em i.e.,} $a_{\sigma} \in
C_K^{\{\dx\}} = C_k^{\{\dx\}} = \CX(t)$. 
Next, a calculation shows that $\dx(\sigma(\log x) - \log x) = 0 = \dt(\sigma(\log x) - \log x)$ so
we have that $\sigma(\log x) = \log x + c_\sigma$ for some $c_\sigma \in \CX$.  Finally, a
calculation shows that
\[
0 = \dt(\sigma(x^t)) - \sigma(\dt (x^t)) = (\dt a_\sigma -a_\sigma c_\sigma) x^t
\]
so we have that 
\begin{eqnarray} 
\dt(\frac{\dt a_\sigma}{a_\sigma} )& =& 0  \ . \label{eqn3}
\end{eqnarray}
Conversely, one can show that for any $a$ such that $\dx a = 0$ and equation (\ref{eqn3}) holds, the
map defined by $x^t \mapsto ax^t, \ \log x \mapsto \log x + \frac{\dt a}{a}$ defines a differential
$k$-automorphism of $K$ so we have
\[
\Gal(K/k) = \{a \in C_K^{(\frac{\d}{\d t})} = C_k^{(\frac{\d}{\d t})} \ | \ a\neq 0 \mbox{ and } \dt(\frac{\dt a}{a})
= 0 \}
\] 
\hfill \QED}
\end{example} This example illustrates two facts. The first is that the Galois group of a
parameterized linear differential equation is a group of $n\times n$ matrices (here $n=1$) whose
entries are functions of the parameters (in this case, $t$) satisfying certain differential
equations; such a group
is called a linear differential algebraic group (see Definition~\ref{LDAG} below). In general, the
Galois group of a parameterized linear differential equation will be such a group.\\
[0.1in]
The second fact is that in this example $\Gal(K/k)$ does not contain enough elements to give a
Galois correspondence. Expressing an element of $\CX(t)$ as $a = a_0\prod (t-b_i)^{n_i}, \ a_0, b_i
\in \CX, n_i \in \ZX$, one can show that if $a \in \Gal(K/k)$ then $a \in \CX$, that is $\Gal(K/k) =
\CX^*$. If $\sigma \in \Gal(K/k)$ and $\sigma(x^t) = ax^t$ with $a \in \CX$, then 
\[
\sigma(\log x) = \sigma(\frac{\dt x^t}{x^t}) = \frac{\dt(ax^t)}{ax^t} = \log x \ .
\]
Therefore $\log x$ is fixed by the Galois group and so there cannot be a Galois 
correspondence. The problem is that we do not have an element $a \in \CX(t)$ such that
$\dt(\frac{\dt a}{a}) = 0$ and $\dt a \neq 0$.\\
[0.1in] In the Picard-Vessiot Theory, one avoids a similar problem by insisting that the constant
subfield is large enough, {\em i.e.,} algebraically closed. This insures that any consistent set of
polynomial equations with constant coefficients will have a solution in the field. In the
Parameterized Picard-Vessiot Theory that we will develop, we will 
need to insure that any consistent system of differential equations (with respect to the parametric
variables) has a solution. This motivates the following definition.\\
[0.1in]
Let $k$ be a $\Delta$-field with derivations $\Delta = \{\d_1, \ldots , \d_m\}$. The {\em
$\Delta$-ring $k\{y_1, \ldots , y_n\}_\Delta$ of differential polynomials in $n$ variables over $k$} is the
usual polynomial ring in the infinite set of variables 
\[
\{\d_1^{n_1}\d_2^{n_2} \cdots \d_m^{n_m}y_j\}_{j=1, \ldots ,n}^{n_i \in \NX}
\]
with derivations extending those in $\Delta$ on $k$ and defined by 
$\d_i(\d_1^{n_1} \cdots \d_i^{n_i}\cdots\d_m^{n_m}y_j) = \d_1^{n_1} \cdots
\d_i^{n_i+1}\cdots\d_m^{n_m}y_j$.
\begin{defin}\label{diffclosed} We say that  a $\Delta$-field $k$ is {\em differentially closed} if
for any $n$ and any set $\{P_1(y_1, \ldots, y_n), \ldots, P_r(y_1, \ldots ,y_n), Q((y_1, \ldots
,y_n)\} \subset k\{y_1, \ldots , y_n\}_\Delta$, if the system 
\[
\{P_1(y_1, \ldots, y_n)=0, \ldots , 
P_r(y_1, \ldots ,y_n)=0, Q(y_1, \ldots ,y_n)\neq 0\}
\]
has a solution in some $\Delta$-field $K$ containing $k$, then it has a solution in $k$
\end{defin} This notion was introduced by A. Robinson \cite{robinson59} and extensively developed by
L. Blum \cite{blum} (in the ordinary differential case) and E.R. Kolchin \cite{constrained} (who
referred to these as constrainedly closed differential fields).  More recent discussions can be found
in
\cite{marker2000} and \cite{mcgrail}.  A fundamental fact is that any $\Delta$-field $k$ is
contained in a differentially closed differential field. In fact, for any 
such $k$ there is a differentially closed $\Delta-$field $\bar k$ containing $k$ such that for any
differentially closed $\Delta$-field $K$ containing $k$, there is a differential $k$-isomorphism of
$\bar k$ into $K$.  Differentially closed fields have many of the same properties with
respect to differential fields as algebraically closed fields have with respect to fields but there
are some striking differences.  For example, the differential closure of a field has proper
subfields that are again differentially closed!  For more information, the reader is referred to the
above papers.
\\
[0.2in]
\noindent{\bf Example 3.1(bis)} Let $k$ be a $\Delta = \{\d_x, \d_t\}$-field and let $k_0 =
C_k^{\dx}$.  Assume that $k_0$ is a differentially closed $\d_t$-field and that $k = k_0(x)$
where $\d_xx = 1$ and $\d_t x = 0$. We again consider the differential equation
\[
\dx y = \frac{t}{x}y \ 
\]
and let $K = k(x^t, \log x)$ where $x^t, \log x$ are considered formally as algebraically
independent elements satisfying $\dt (x^t) = (\log x)x^t, \
\dx (x^t) = \frac{t}{x}x^t, \ \dt(\log x) = 0, \dx (\log x) = \frac{1}{x}$. One can show that
$C_K^{\{\dt\}} = k_0$ and that the 
Galois group is again
\[
\Gal(K/k) = \{ a \in k_o^* \ | \ \dt(\frac{\dt(a)}{a}) = 0 \} \ .
\]
Note that $\Gal(K/k)$  contains an element $a$ such that $\dt a \neq 0 $ and $\dt(\frac{\dt(a)}{a})
= 0$.  To see this, 
note that the $\{\d_t\}$-field $k_0(u)$, where $u$ is transcendental over $k_0$ and $\dt u = u$ is
a $\{\d_t\}$-extension of $k_0$ containing such an element ({\em e.g.,} $u$).  The definition of
differentially closed ensures that $k_0$ also contains such an element. This implies that $\log x$
is not left fixed by $\Gal(K/k)$.  In fact, we will  show in Section \ref{multgps} that the
following is a complete list of differential algebraic subgroups of $\PGal(K/k)$ and the corresponding $\Delta$-subfields
of $K$:
\[
\begin{array}{|cc|}\hline {\rm Field} & {\rm Group} \\
\hline 
 k((x^t)^n, \log x), z\in \NX_{>0} & \{a \in k_0^* \ | \ a^n = 1\} = \ZX/n\ZX \\
[.05in] 
k( \log x) & \{ a \in k_0^* \ | \ \dt a = 0 \mbox{ and } \dt(\dt(a)/a) = 0 \} \\
[.05in]
 k & \{  a \in k_0^* \ | \ \dt(\dt(a)/a) = 0 \}\\
\hline
\end{array}
\]
\hfill \QED

\noindent We now turn to stating the  Fundamental Theorem in the Galois Theory of Parameterized Linear Differential Equations.
 We need to give a formal definition of the kinds of groups that can occur and also 
of what takes the place of a Picard-Vessiot extension.  This is done in the next two definitions.

\begin{defin}\label{LDAG} Let $k$ be a differentially closed $\Delta$-differential field.\\
\indent 1.~ A subset $X \subset k^n$ is said to be {\em Kolchin-closed} if there exists a set $\{ 
f_1, \ldots , f_r\}$ of differential polynomials in $n$ variables such that $X = \{ a \in k^n \ | \
f_1(a) = \ldots = f_r(a) = 0\}$.\\
\indent 2.~ A subgroup  $G \subset \GL_n(k) \subset k^{n^2}$ is a {\em linear differential algebraic
group} if $G = X \cap \GL_n(k)$ for some Kolchin-closed subset of $k^{n^2}$.
\end{defin}
In the previous example, the Galois group was exhibited as a linear differential algebraic subgroup
of $\GL_1(k_0)$. For any linear algebraic group $G$, the group $G(k)$  is a linear differential
algebraic group. Furthermore, the group $G(C^\Delta_k)$ of constant points of $G$ is also a linear
differential algebraic group since it is defined by the (algebraic) equations defining $G$ as well
as the (differential) equations stating that the entries of the matrices are constants.  We will
give more examples in the next section\\
[0.1in]
In the next definition, we will use the following conventions. If $F$ is a $\Delta = \{\d_0, \d_1,
\ldots , \d_m\}$-field, we denote by $C_F^0$ the $\d_0$ constants of $F$, that is, $C_F^0 =
C_F^{\{\d_0\}} = \{c \in F \ | \ \d_0c = 0\}$. One sees that $C_F^0$ is a $\Pi = \{
\d_1, \ldots , \d_m\}$-field. We will use the notation $k\langle z_1, \ldots , z_r \rangle_\Delta$
to denote a $\Delta$-field containing $k$ and elements $z_1, \ldots ,z_r$ such that no proper
$\Delta$-field has this property, {\em i.e.,}$k\langle z_1, \ldots , z_r \rangle_\Delta$ is the field generated over $k$ by $z_1, \ldots ,z_n$
and their higher derivatives. 

\begin{defin} Let $k$ be a $\Delta = \{\d_0, \d_1, \ldots , \d_m\}$-field and let 
\[
\d_0Y = A Y
\]
be a differential equation with $A \in \M_n(k)$. 
\begin{enumerate}
\item A $\Delta$-extension $K$ of $k$ is a 
{\em Parameterized Picard-Vessiot extension of $k$ (or, more compactly, a PPV-extension of $k$)} if
$K = k\langle z_{1,1}, \ldots ,z_{n,n}\rangle_\Delta$ where 
\begin{enumerate}
\item the $z_{i,j}$ are entries of a matrix $Z \in \GL_n(K)$ satisfying $\d_0Z = AZ$, and
\item $C^0_K = C^0_k$, {\em i.e.,} the $\d_0$-constants of $K$ coincide with the $\d_0$-constants of
$k$.
\end{enumerate}
\item
The group $\PGal(K/k) = \{\sigma : K\rightarrow K \ | \ \sigma \mbox{ is a }
k\mbox{-automorphism such that } \sigma \d = \d \sigma \ \forall  \d \in \Delta \}$ is called
the {\em Parameterized Picard-Vessiot Group (PPV-Group)} associated with the PPV-extension $K$ of
$k$.\end{enumerate}
\end{defin}
We note that if $K$ is a PPV-extension of $k$ and $Z$ is as above then for any $\sigma \in
\PGal(K/k)$ one has that $\d_0(\sigma(Z)Z^{-1}) = 0$.  Therefore we can identify each $\sigma \in
\PGal(K/k)$ with a matrix in $GL_n(C^0_k)$. 
We can now state the Fundamental Theorem of Parameterized Picard-Vessiot Extensions
\begin{thm}\label{PPVthm} Let $k$ be a $\Delta = \{\d_0, \d_1, \ldots , \d_m\}$-field and assume
that $C_k^0$ is a differentially closed $\Pi=\{\d_1, \ldots ,\d_m\}$-field. Let
\begin{eqnarray}
 \d_0Y = A Y  \label{fundeqn}
 \end{eqnarray}
be a differential equation with $A \in \M_n(k)$.
\begin{enumerate}
\item There exists a PPV-extension $K$ of $k$ associated with (\ref{fundeqn}) and this is unique up
to differential $k$-isomorphism.
\item The PPV-Group $\PGal(K/k)$ equals $G(C_k^0)$, where $G$ is a linear differential algebraic group defined over $C_k^0$.
\item The  map that sends any $\Delta$-subfield $F, k \subset F \subset K$, to the group
$\PGal(K/F)$ is a bijection between differential subfields of $K$ containing $k$ and Kolchin-closed
subgroups of $\PGal(K/k)$. Its inverse is given by the map that sends a Kolchin-closed group $H$ to
the field $\{z \in K \ | \ \sigma(z) = z \mbox{ for all } \sigma \in H\}$.
\item A Kolchin-closed subgroup $H$ of $\PGal(K/k)$ is a normal subgroup of $\PGal(K/k)$ if and
only if the field $K^H$ is left set-wise invariant by   $\PGal(K/k)$. If this is the case, the 
map $\PGal(K/k) \rightarrow  \PGal(K^H/k)$ is surjective with kernel $H$ and
$K^H$ is a PPV-extension of $k$ with PPV-group isomorphic to $\PGal(K/k)/H$. Conversely, if $F$ is
a differential subfield of $K$ containing $k$ and $F$ is a PPV-extension of $k$, then $\PGal(F/K)$
is a normal subgroup of $\PGal(K/k)$. 
\end{enumerate}
\end{thm} 
The proof of this result is virtually the same as for the corresponding result of Picard-Vessiot
Theory. We give the details in Appendices 9.1-9.4.\\[0.1in]
We
will give two simple applications of this theorem.
 For the first, let $K$ be a PPV-extension of $k$ corresponding to the equation $\d_0Y=AY$ and let $K =k\langle
Z\rangle_\Delta$, where $Z \in \GL_n(K)$ and $\d_0 Z= AZ$. We now consider the field  $\KPV = k(Z) \subset K$.  
Note that $\KPV$ is not necessarily
a $\Delta$-field but it is a $\{\d_0\}$-field.  One can easily see that it is a PV-extension for
the equation $\d_0Y - AY$ and that the PPV-group leaves it invariant and acts as
$\{\d_0\}$-automorphisms. We therefore have an injective homomorphism of $\PGal(K/k) \rightarrow
{\rm Gal}_{\{\d_0\}}(\KPV/k)$, defined by restriction.   We then have the following result
\begin{prop}\label{PPVvsPV} Let $k, C_k^0, K$ and $\KPV$ be as above. Then:
\begin{enumerate}
\item When considered as ordinary $\{\d_0\}$-fields, $\KPV$ is a PV-extension of $k$ with
algebraically closed field $C_k^0$ of $\d_0$-constants.
\item If ${\rm Gal}_{\{\d_0\}}(\KPV/k)\subset \GL_n(C_k^0)$ is the Galois group of the  ordinary differential field
$\KPV$ over $k$, then
the Zariski closure of the Galois group $\PGal(K/k)$ in $\GL_n(C_k^0)$ equals ${\rm Gal}_{\{\d_0\}}(\KPV/k)$.
\end{enumerate}
\end{prop}
\begin{proof} Since a differentially closed field is algebraically closed, we have already justified
the first statement. Clearly, $\PGal(K/k)
\subset{\rm Gal}_{\d_0}(K/k)$. Since $\PGal(K/k)$ and ${\rm Gal}_{\{\d_0\}}(\KPV/k)$ have the same fixed field $k$, the second
statement follows.
\end{proof}

\begin{remark}\label{KPVremark}{\em One sees that the field $\KPV$ is independent of the particular invertible solution $Z$ of $\d_0Y = AY$ used to generate $K$ (although the Galois groups are only determined 
up to conjugacy). On the other hand,  $\KPV$  does depend on the the equation $\d_0Y = AY$ and not just on the field $K$, that is if $K$ is a PPV-extension of 
$k$ for two different equations $\d_0Y = A_1Y$ and $\d_0Y = A_2Y$ with solutions $Z_1$ and $Z_2$ respectively, the fields $K^{\rm PV}_{A_1}$ and $K^{\rm PV}_{A_2}$ (and their respective PV-groups)
may be very different.  We will give an example of this in Remark~\ref{gamma}.}

\end{remark}
Our second application is to characterize those equations $\d_0Y = AY$ whose PPV-groups are the set
of $\Delta$-constant points of a linear algebraic group. We first make the following definition.
\begin{defin} Let $k$ be a $\Delta$-differential field and let 
$A \in \M_n(k)$.  We say that $\d_0Y=AY$ is {\em completely integrable} if there exist $A_i \in
\M_n(k), \ i = 0, \ldots, n$ with $A_0 = A$ such that 
\begin{eqnarray*}
\d_j A_i - \d_i A_j = A_jA_i - A_iA_j  \ \mbox{ for all $i,j = 0, \ldots$ n.} 
\end{eqnarray*}
\end{defin} The latter conditions on the $A_i$ are the usual integrability conditions and motivate
the nomenclature.
\begin{prop} \label{propint} Let $k$ be a $\Delta$-differential field and assume that 
$k_0$ is a $\Pi$-differentially closed $\Pi$-field.  Let $A \in \M_n(k)$ and let $K$ 
be a PPV-extension of $k$ for $\d_0Y = AY$. Finally, let $C = 
C_{k}^{\Delta}$.
\begin{enumerate}
\item There exists a linear algebraic group $G$ defined over $C$ such that $\PGal(K/k)$ is
conjugate to  $G(C)$ if and only if $\d_0Y = AY$ is completely integrable. If this is the case,
then $K$ is a PV-extension of $k$ corresponding to this integrable system.
\item If $A \in \M_n(C_k^{\Pi})$, then $\PGal(K/k) \mbox{ is conjugate to } G(C)$ for some
linear algebraic group defined over $C$.
\end{enumerate}
\end{prop}
\begin{proof} 1.~Let $K = k\langle Z \rangle_\Delta$ where $Z \in \GL_n(K)$ satisfies $\d_0Z = AZ$.  
If the PPV-group  is as described, then  there exists a $B \in \GL_n(C_k^0)$ such that $B\PGal(K/k)B^{-1} = G(C)$, $G$ 
an algebraic subgroup of $\GL_n(C_k^0)$, defined over $C$. Set $W = ZB^{-1}$. One sees that $\d_0W = AW $ and $K= k\langle W \rangle_\Delta$.
A simple calculation shows that  
for any $i = 0, \ldots , n$, $\d_iW\cdot W^{-1}$ is left fixed by all elements of the PPV-group.
Therefore $\d_i W = A_i W$ for some $A_i \in \M_n(k)$. Since the $\d_i$ commute, one sees that the
$A_i$ satisfy the integrability conditions.  \\
[0.1in]
Now assume that there exist $A_i \in \M_n(k)$ as above satisfying the integrability conditions.
Let $K$ be a PV-extension of $k$ for the corresponding integrable system. From
Lemma~\ref{PV=PPV} in the Appendix, we know that $K$ is also a PPV-extension of $k$ for $\d_0Y
= AY$. Let $\sigma
\in
\PGal(K/k)$ and let $\sigma(Z) = Z D$ for some $D \in \GL_n(C_k^0)$.  Since $\d_iZ\cdot Z^{-1}  = 
A_i\in
\GL_n(k)$, we have that $\sigma(\d_iZ\cdot Z^{-1}) = \d_iZ\cdot Z^{-1}$. A calculation then shows
that $\d_i(D) = 0$.  Therefore $D \in \GL_n(C_K^\Delta)$. We now need to show that $C_K^\Delta =
C_k^\Delta$.  This is clear since $C_K^\Delta  \subset C_K^0 = C_k^0$. The final claim of part 1. is
now clear.
\\[0.1in]
 2.~Under the assumptions, the matrices $A_0 = A, A_1 = 0, \ldots , A_n = 0$ satisfy the
integrability conditions, so the  conclusion follows from Part 1. above.
\end{proof}
If $A$ has entries that are analytic functions of $x, t_1, \ldots ,t_m$, the fact that $\PGal(K/k) =
G(C)$ for some linear algebraic group does not imply that, for  some open set
of values $\vta =(\tau_1,
\ldots ,\tau_m)$ of $(t_1, \ldots ,t_m)$, the Galois group $G_\vta$ of the ordinary differential
equation $\d_xY = A(x,\tau_1, \ldots ,\tau_m)Y$ is  independent of the choice of $\vta$. We shall 
see in Section~\ref{isomonsec} that for equations with regular singular points we do have a constant
Galois group (on some open set of parameters)  if the PPV-group is $G(C)$ for some linear algebraic
group but the following shows that this is not true in general.\\

\begin{example}\label{isogalois}{\em Let $\Pi = \{\d_1 = \frac{\d}{\d t_1}, \d_2 = 
\frac{\d}{\d t_2}\}$ and $k_0$ be a differentially closed $\Pi$-field containing $\CX$.  Let $k =
k_0(x)$ be a $\Delta = \{\d_0 = \frac{\d}{\d x}, \d_1, \d_2\}$-field where $\d_0|_{k_0} = 0, \d_0(x)
= 1,$ and $\d_1, \d_2$ extend the derivations on $k_0$ and satisfy $\d_1(x) = \d_2(x) = 0$. 
The equation
\begin{eqnarray*}
\frac{\d Y}{\d x}& = A(x,t_1, t_2)Y &=\left(
\begin{array}{cc}t_1 & 0 \\
0 & t_2
\end{array}\right)Y
\end{eqnarray*}
has solution
\begin{eqnarray*}
Y &= &\left(
\begin{array}{cc}e^{t_1x} & 0 \\
0 & e^{t_2x}
\end{array}\right)
\end{eqnarray*}
One easily checks that 
\[
A_1 = \frac{\d Y}{\d t_1} Y^{-1} \in \gl_2(k) \ \ \mbox{ and }\ \ A_2 = 
\frac{\d Y}{\d t_2} Y^{-1}
\in
\gl_2(k) 
\]
so the Galois group associated to this equation is conjugate to $G(C)$ for some linear algebraic
group $G$ (in fact
$G(C) = C^* \times C^*$).  Nonetheless, for fixed values $\vta = (\tau_1, \tau_2) \in C^2$, the
Galois group of $\d_0Y = A(x,\tau_1,\tau_2)Y$ is $G(C)$ if and only if the only rational numbers
$r_1, r_2$ satisfying $r_1\tau_1 + r_2\tau_2 = 0$ are $r_1 = r_2 = 0$. \hfill \QED}
\end{example}
For more information on how a differential Galois group can vary in a family of linear differential
equations see  \cite{andre2004}, \cite{berken02}, \cite{berken04}, \cite{Hrushovski},
 and
\cite{singer_moduli}.\\
[0.1in]
We end this section with a result concerning solving parameterized linear differential eqautions in
``finite terms''. The statement of the result is the same {\em mutatis mutandi} as the corresponding
result
in the usual Picard-Vessiot Theory ({\em cf.,} \cite{PuSi2003}, Ch. 1.5) and will be proved in the
Appendix.
\begin{defin}\label{lioudef} Let $k$ be a $\Delta = \{\d_0, \ldots \d_m\}$-field.  We say that a
$\Delta$-field $L$ is a {\em parameterized liouvillian extension} of $k$ if $C_L^0 = C_k^0$ and
there exist a tower of $\Delta-$fields $k = L_0 \subset L_1 \subset \ldots \subset L_r$ such that
$L_i = L_{i-1}\langle t_i \rangle_\Delta$ for $i = 1 \ldots r$, where either
\begin{enumerate}
\item $\d_0 t_i \in L_{i-1}$, that is $t_i$ is a {\em parameterized integral (of an element of
$L_{i-1}$)}, or
\item $t_i \neq 0 \mbox{ and } \d_0 t_i / t_i \in L_{i-1}$, that is $t_i$ is a {\em parameterized exponential
(of an integral of an element in $L_{i-1}$)}, or
\item $t_i$ is algebraic over $L_{i-1}$.
\end{enumerate}
\end{defin}
One then has the following result
\begin{thm}\label{liouvillian} Let $k$ be a $\Delta$-field and assume that $C^0_k$ is a
differentially closed $\Pi = \{\d_1, \ldots \d_m\}$-field. Let $K$ be a PPV-extension of $k$
with PPV-group $G$. The following are equivalent
\begin{enumerate}
\item $G$ contains a solvable subgroup $H$ of finite index.
\item $K$ is a parameterized liouvillian extension of $k$.
\item $K$ is contained in a parameterized liouvillian extension of $k$.
\end{enumerate}
\end{thm}

\section{Linear Differential Algebraic Groups}\label{multgps} In this section we review some known facts concerning
linear differential algebraic groups and give some examples of these groups.  The theory of linear
differential algebraic groups was initiated by P. Cassidy in \cite{cassidy1} and further developed
in \cite{cassidy2,cassidy3,cassidy4,cassidy5,cassidy6}.  The topic has 
also been addressed in \cite{buium_groups}, \cite{kolchin_groups}, \cite{pillay1996}, 
\cite{pillay97},
\cite{pillay2000},
\cite{sit74}, and \cite{sit75}.  For a general overview see \cite{buium_cassidy}.\\
[0.1in]
Let $k_0$ be a differentially closed $\Delta'= \{\d_1, \ldots , \d_m\}$-field and let $C =
C_{k_0}^{\Delta'}$. 
As we have already defined,
 a linear differential algebraic group is a Kolchin-closed subgroup of $\GL_n(k_0)$.  Although the 
definition is a natural generalization of the definition of a linear algebraic group there are many
points at which the theories diverge.  The first is that an affine differential algebraic group (a
Kolchin-closed subset $X$ of $k_0^m$ with group operations defined by everywhere defined rational
differential functions) need not be a linear differential algebraic group although affine
differential algebraic groups whose group laws are given by differential polynomial maps are linear
differential algebraic groups
\cite{cassidy1}. Other distinguishing phenomena will emerge as we examine some examples.\\[0.2in]
{\bf Differential Algebraic Subgroups of $\Ga^n$.} The group $\Ga = (k_0,+)$ is naturally
isomorphic to $\{\left(
\begin{array}{cc} 1 & a \\
0& 1 
\end{array}\right) \ | \ a \in k_0 \}$ and, as such, has the structure of a linear differential
algebraic group. Nonetheless we will continue to identify this group with $k_0$. The set $\Ga^n =
(k_0^n, +)$ can also be seen to be a linear differential algebraic group. In (\cite{cassidy1}, Lemma
11), Cassidy shows that a subgroup $H$ of $\Ga^n$ is a linear differential algebraic group if and
only
if $H$ is the set of zeros of a set of linear homogeneous differential polynomials in $k_0\{y_1, 
\ldots , y_n\}$.  In particular, when $m = n = 1, \Delta' =
\{\d\}$, 
the linear differential algebraic subgroups of $\Ga$ are all of the form 
\[
\Ga^L(k_0) = \{ a \in \Ga(k_0) \ | \ L(a) = 0 \}  
\]
where $L$ is a linear differential operator ({\em i.e.,} an element of the ring $k_0[\d]$ whose
multiplication is given by $\d \cdot a = a\d + \d(a)$). 
The lattice structure of these subgroups is given by
\[
\Ga^{L_1}(k_0) \subset  \Ga^{L_2}(k_0) \Leftrightarrow L_2 = L_3L_1\mbox{ for some } L_3 \in k_0[\d]
\
. 
\]
\vspace{.1in}

\noindent {\bf Differential Algebraic Subgroups of $\Gm^n$} These have been classified
by Cassidy (\cite{cassidy1}, Ch.IV).We shall restrict ourselves to the case $n = m = 1, \Delta' =
\{\d\}$, that is, differential algebraic subgroups of $\Gm(k_0) = \GL_1(k_0) = k_0^*$.
 Any such group is either
\begin{enumerate}
\item finite and cyclic, or
\item $\Gm^L = \{a \in \Gm(k_0) \ | \ L(\frac{\d a}{a}) = 0\}$ for some $L \in k_0[\d].\}$
\end{enumerate} For example, if $L = \d$, the group 
\[
\Ga^{\d}(k) = \{ a \in k_0^* \ | \ \d(\frac{\d a}{a}) = 0 \}
\]
is the PPV-group of the parameterized linear differential equation $\d_x y = \frac{t}{x}y$ where $\d
= \d_t$. Notice that the only proper differential algebraic subgroup of $\{ a \in k_0 \ | \d a = 0\}$  is
$\{0\}$.  Therefore the only proper differential algebraic subgroups of $\Ga^{\d}$ are either the finite
cyclic groups, or $\Gm(C)$.  This justifies the left column in the table given in Example 3.1(bis).
 The right column follows by calculation.\\
[0.1in]
The proof that the groups of $1.$ and $2.$ are the only possibilities proceeds in two steps. The
first is to show that if the group is not connected (in the Kolchin topology where closed sets are
Kolchin-closed sets), it must be finite (and therefore cyclic). The second step involves the {\em
logarithmic derivative map} $\ld: \Gm(k_0) \rightarrow \Ga(k_0)$ defined by 
\[
\ld(a) = \frac{\d a}{a} \ . 
\]
This map is a differential rational map ({\em i.e.,} the quotient of differential polynomials) and
is a homomorphism. 
Furthermore, it can be shown that the following sequence is exact
\[
(1)\rightarrow \Gm(C) \rightarrow \Gm(k_0) \stackrel{a\mapsto \frac{\d a}{a}} {\rightarrow}\Ga(k_0)
\rightarrow (0) \ .
\]
The result then follows from the classification of differential subgroups of $\Ga(k_0)$. Note that in
the usual theory of linear algebraic groups, there are no nontrivial rational homomorphisms from
$\Gm$ to $\Ga$.\\[0.2in]
{\bf Semisimple Differential Algebraic Groups} These groups have been classified by Cassidy
in \cite{cassidy6}.  Buium \cite{buium_groups} and Pillay \cite{pillay97} have given simplified
proofs in the ordinary case ({\em i.e.,} $m = 1$). Buium's proof is geometric using the notion of jet
groups and Pillay's proof is model theoretic and assumes from the start that the groups are finite
dimensional (of finite Morely rank).\\
[0.1in] We say that a connected differential algebraic group is semisimple if it has no nontrivial
Kolchin-connected, commutative subgroups. Let us start by considering semisimple differential
algebraic subgroups $G$ of $\SL_2(k_0)$.   Let $H$ be the Zariski-closure of such a group. If $H
\neq \SL_2(k_0)$, then $H$ is solvable ({\em cf.,}
\cite{PuSi2003}, p.127) and so the same would be true of $G$.  Therefore $G$ must be Zariski-dense
in $\SL_2(k_0)$.  In \cite{cassidy1}, Cassidy classified the Zariski-dense differential algebraic subgroups
of $\SL_n(k_0)$.  Let $\DX$ be the $k_0$-vector space of derivations spanned by $\Delta'$. 
\begin{prop}\label{sln} Let $G$ be a proper Zariski-dense differential algebraic subgroup of $\SL_n(k_0)$.  Then there
exists a finite set $\Delta_1 \subset \DX$ of commuting derivations such that $G$ is conjugate to
$\SL_n(C_{k_0}^{\Delta_1})$, the $\Delta_1$-constant points of $\SL_n(k_0)$. 
\end{prop} Note that in the ordinary case  $m=1$, we can restate this more simply: {\em A proper
Zariski-dense subgroup of $\SL_n(k_0)$ is conjugate to $\SL_n(C)$}. A complete classification of
differential subgroups of $\SL_2$ is given in \cite{sit75}. The complete classification of
semisimple differential algebraic groups  is given by the following result. By a Chevalley group, we
mean a connected simple $\QX$-group containing a maximal torus diagonalizable over $\QX$. 
\begin{prop} Let $G$ be a Kolchin-connected semisimple linear differential algebraic group. Then
there exist finite subsets of commuting derivations $\Delta_1, \ldots , \Delta_r$ of $\DX$, 
Chevalley groups $H_1, \ldots, H_r$ and a differential isogeny $\sigma:
H_1(C_{k_0}^{\Delta_1})\times 
\cdots \times H_r(C_{k_0}^{\Delta_r}) \rightarrow G$ \footnote{One need not assume that $G$ is
linear since Pillay \cite{pillay1996} showed that a semisimple differential algebraic group is
differentially isomorphic to a linear differential algebraic group.}.
\end{prop}
 
\section{Isomonodromic  Families} \label{isomonsec} In this section we shall describe how
isomonodromic families of linear differential equations fit into this theory of parameterized linear
differential equations. We begin with some definitions and follow the exposition of Sibuya
\cite{sibuya}, Appendix 5\footnote{The presentation clearly could be cast in the language of vector
bundles (see
\cite{berken02}, \cite{berken04},
\cite{mal80a}, \cite{mal80b}) but the approach presented here is more in the spirit
of the rest of the paper.}. Let $\calD$ be an open subset of the Riemann sphere (for simplicity, we
assume that the point at infinity is not in $\calD$) and let $\calD(\vta,\vr) =
\prod_{h=1}^pD(\tau_h,\rho_h)$  where $\vr = (\rho_1,
\ldots , \rho_p)$ is a $p$-tuple of positive numbers, $\vta = (\tau_1, \ldots , \tau_p) \in
\CX^p$ and  $D(\tau_h,\rho_h)$ is the open disk in $\CX$ of radius $\rho_h$ centered at the point
$\tau_h$.  We denote by $\calO(\calD\times\calD(\vta,\vr))$ 
the ring of functions $f(x,\vt)$ holomorphic on $\calD\times\calD(\vta,\vr)$. Let $A(x,\vt) \in
\M_n(
\calO(\calD\times\calD(\vta,\vr))$ and consider the differential equation
\begin{eqnarray}
\frac{\d W}{\d x} = A(x,\vt)W \label{isoeqn}
\end{eqnarray}
\begin{defin} A {\em system of fundamental solutions} of (\ref{isoeqn})  is a collection of pairs
$\{D(x_j,\vr_j),W(x,\vt)\}$ such that
\begin{enumerate}
\item the disks $D(x_j,\vr_j)$ cover $\calD$ and
\item for each $\vt \in \calD(\vta,\vr)$ the $W_j(x,\vt) \in
\GL_n(\calO(D(x_j,\vr_j)\times\calD(\vta,\vr))$ are solutions of (\ref{isoeqn}). 
\end{enumerate}

We define $C_{i,j}(\vt) = W_i(x,\vt)^{-1} W_j(x,\vt)$ whenever $D(x_i,\vr_i)\cap
D(x_j,\vr_j)\neq\emptyset$ and refer to these as the {\em connection matrices} of the system of
fundamental solutions.
\end{defin}
\begin{defin} The differential equation (\ref{isoeqn}) is {\em isomonodromic on
$\calD\times\calD(\vta,\vr)$} if there exists a system $\{D(x_j,\vr_j),W_j(x,\vt)\}$ of fundamental
solutions such that the connection matrices $C_{i,j}(\vt)$ are independent of $t$. 
\end{defin}
We note that for a differential equation that is isomonodromic in the above sense, the monodromy
around any path is independent of $\vt$ as well.  To see this let $\gamma$ be a path in $\calD$
beginning and ending at $x_0$ and let $D(x_1,\vr_1), \ldots , D(x_s,\vr_s),  D(x_1,\vr_1)$ be a
sequence of disks covering the path so that $D(x_i,\vr_i)\cap D(x_{i+1},\vr_{i+1}) \neq
\emptyset$ and $x_0 \in D(x_1,\vr_1)$. If one continues $W_1(x_1,\vt)$ around the path then the
resulting matrix $\tilde{W} =W_1(x_1,\vt) C_{1,s}C_{s,s-1}\cdots C_{2,1}$. By assumption, the
monodromy matrix $C_{1,s}C_{s,s-1}\cdots C_{2,1}$ is independent of $\vt$. \\
[0.1in]
For equations with regular singular points, the monodromy group is Zariski dense in the PV-group.
The above  comments therefore imply that for an isomonodromic family, there is a nonempty
open set of paramters such that for these values the PV-group is constant as the parameters vary in
this set. As Sibuya points out, it is not clear if the reverse implication is true in general
(although it is true in special cases \cite{boli97}) . \\
[0.1in]
With these definitions, Sibuya shows (\cite{sibuya}, Theorem A.5.2.3)
\begin{prop}\label{sibuya} The differential equation (\ref{isoeqn}) is isomonodromic on
$\calD\times\calD(\vta,\vr)$ if and only if there exist $p$ matrices $B_h(x,\vt) \in \gl_n(
\calO(\calD\times\calD(\vta,\vr))), \ h = 1,\ldots ,p$ such that the system
\begin{eqnarray}
\frac{\d W}{\d x} & = & A(x,\vt) W \notag \\
\frac{\d W}{\d t_h} & = & B_h(z,\vt)W \ \ \ (h = 1, \ldots ,p) \label{complete}
\end{eqnarray}
is completely integrable.
\end{prop} Some authors use the existence of matrices $B_i$ as in Proposition~\ref{sibuya} as the
definition of isomonodromic ({\em cf.,}\cite{mal80a}). Sibuya goes on to note that if $A(x,\vt)$ is
rational in $x$ and if the differential equation has only regular singular points, then the
$B_h(x,\vt)$ will be rational in $x$ as well (without the assumption of regular singular points 
one cannot conclude that the $B_i$ will be rational in $x$.) This observation leads to the next
proposition. 
\\
[0.1in]
For any open set $\calU \subset \CX^p$, let $\calM(\calU)$ be the field of functions meromorphic on
$\calU$. Note that $\calM(\calU)$ is a $\Pi = \{\frac{\d}{\dt_1},\ldots
,\frac{\d}{\dt_p}\}$-field. If $\calU'
\subset
\calU$ then there is a natural injection of ${\rm res}_{\calU,
\calU'}:\calM(\calU)\rightarrow \calM(\calU')$. We shall need the following result 
of Seidenberg
\cite{sei58, sei69}: {\em Let $\calU$ be an open subset of $\CX^p$ and let $F$ be a
$\Pi$-subfield of $\calM(\calU)$ containing $\CX$.  If $E$ is $\Pi$-field containing $F$ and
finitely generated (as a $\Pi$-field) over $\QX$, then there exists a nonempty open set $\calU'
\subset \calU$ and an isomorphism $\phi: E \rightarrow \calM(\calU)$ such that $\phi|_F = {\rm
res}_{\calU,
\calU'}$}.\\
[0.1in]
Let $A(x,\vt)$ be as above, assume the entries of $A$ are rational in $x$ and let $F$ be the
$\Pi$-field generated by the coefficients of powers of $x$ that appear in the entries of $A$.
Let $k_0$ be the differential closure of $F$. We consider $k = k_0(x)$ to be a $\Delta =
\{\frac{\d \ }{\d x}, \frac{\d \ }{\d t_1},\ldots ,\frac{\d \ }{\d t_p}\}$-field in the obvious
way.
Given open subsets $U_1 \subset U_2$ of the Riemann Sphere, we say that $u_1$ is a {\em punctured
subset of $U_2$} if there exist a finite number of disjoint closed disks $D_1, \ldots D_r \subset U_2$ such
that $U_1 = U_2\backslash (\cup_{i=1}^r D_i)$.
\begin{prop}\label{isom3}Let $A(x,\vt), k_0$ and $k$  be as above.  Assume that the differential
equation
\begin{eqnarray}
\frac{\d W}{\d x} = A(x,\vt)W \label{isomeqn2}
\end{eqnarray}
has only regular singular points. Then this equation is isomonodromic on $\calD'\times\calU$, for
some nonempty, open $\calU \subset \calD(\vta,\vr)$ and $\calD'$ a punctured subset of $D$ if and
only if the PPV-group of this
equation over $k$ is of the form $G(\CX)$ for some linear algebraic group $G$ defined over $\CX$.
In this case,  the monodromy group of $(\ref{isomeqn2})$ is independent of
$\Vec t\in U$.
\end{prop}
\begin{proof} Assume that (\ref{isomeqn2}) is isomonodromic. Proposition~\ref{sibuya} and the
comments after it ensure that we can complete (\ref{isomeqn2}) to a completely integrable system
(\ref{complete}) where the $B_i(x,\vt)$ are rational in $x$.  The fact that this is a completely
integrable system is equivalent to the coefficients of the powers of $x$ appearing in the entries of
the $B_i$ satisfying a system $\calS$ of $\Pi$-differential equations with coefficients in
$k_0$.  Since this system has a solution and $k_0$ is differentially closed, the system must have a
solution in $k_0$. Therefore we may assume that the $B_i \in \gl_n(k)$. An application of
Proposition~\ref{propint}.1 yields the conclusion.\\
[0.1in] Now assume  that the PPV-group is $G(\CX)$ for some linear algebraic group $G$.
Proposition~\ref{propint}.1 implies that we can complete (\ref{isomeqn2}) to a completely integrable
system (\ref{complete}) where the $B_i(x,\vt)$ are in $\gl_n(k)$.  Let $E$ be the $\Pi$-field
generated by the coefficients of powers of $x$ appearing in the entries of $A$ and the $B_i$. By the
result of Seidenberg referred to above, there is a nonempty open set $\calU \subset
\calD(\vta,\vr)$ such that  these coefficients can be assumed to be analytic on $\calU$.The
matrices $B_i$ have entries that are rational in $x$ and so may have poles (depending on $\vt$) in
$D$.
By shrinking $\calU$ if necessary and replacing $D$ with a punctured subset $D'$ of $D$, we can
assume that $A$ and the $B_i$ have entries that are holomorphic in $D' \times \calU$. We now
 apply
Proposition~\ref{sibuya} to reach the conclusion.
\end{proof}

\section{Second Order Systems}\label{secondsec}In this section we will apply the results of the previous four
sections to give a classification of parameterized second order systems of linear differential
equations.   We will first consider the case of second order parameterized linear equations
depending on only one parameter.  
\begin{prop} \label{prop2x2} Let $k$ be a $\Delta=\{\d_0,  \d_1\}$-field, assume that $k_0 =C^0_k$
is a differentially closed $\Pi=\{\d_1\}$-field and let $C = C^\Delta_k$. Let $A\in\sl_2(k)$ and
let $K$ be the PPV-extension corresponding to the differential equation 
\begin{eqnarray}
\d_0Y = AY \ . \label{2x2}
\end{eqnarray}
Then, either
\begin{enumerate}
\item $\PGal(K/k)$ equals   $\SL_2(k_0)$, or
\item  $\PGal(K/k)$ contains a solvable subgroup of finite index and $K$ is a parameterized
liouvillian extension of $k$, or
\item  $\PGal(K/k)$ is conjugate to $\SL_2(C)$ and there exist 
 $B_1 \in \SL_2(k)$ such that the system
\begin{eqnarray*}
\d_0 Y &= & AY \\
\d_1 Y & = & B_1 Y
\end{eqnarray*}
 is an integrable system.
 \end{enumerate}
\end{prop} 
\begin{proof} Let $Z \in \GL_2(k)$ be a fundamental solution matrix of (\ref{2x2}) and let $z =
\det Z$. We have that $\d_0 z = (\mbox{\rm trace}A)z$ (\cite{PuSi2003}, Exercise 1.14.5), so $z \in
k_0$. For any $\sigma
\in
\PGal(K/k)$, $z = \sigma (z) = \det \sigma \cdot z$ so $\det \sigma = 1$.  Therefore, $\PGal
\subset \SL_2(k_0)$. Let $G$ be the Zariski-closure of $\PGal(K/k)$. If $G \neq \SL_2(k_0)$,
then $G$ has a solvable subgroup of finite index and so the same holds for $\PGal(K/k)$. 
Therefore, 2.~holds.
If $ G =\SL_2(k_0)$, then 1.~holds.   If $ G = \SL_2(k) \mbox{ and } G \neq \PGal(K/k)$, then
Proposition~\ref{sln} and the discussion immediately following it imply that there is a $B \in
\SL_2(k_0)$ such that $B\PGal(K/k)B^{-1} = \SL_2(C)$. Proposition~\ref{propint} then implies that
the parameterized equation $\d_0Y = AY$ is completely integrable, yielding conclusion 3.
\end{proof}
If the entries of $A$ are functions of $x$ and $t$, analytic in some domain and rational in $ x$, we
can combine the above proposition with Proposition~\ref{isom3} to yield the next corollary. Let
$\calD$ be an open region on the Riemann Sphere and $D(\tau_0,\rho_0)$ be the open disk of radius 
$\rho_0$ centered at $\tau_0$ in $\CX$.  Let $\calO(\calD \times D(\tau_0,\rho_0)$ be the ring of
functions holomorphic in $\calD \times D(\tau_0,\rho_0)$ and let $A(t,x) \in \sl_2(\calO(\calD \times
D(\tau_0,\rho_0))$ and assume that $A(x,t)$ is rational in $x$. Let 
$\Delta = \{\d_0 = \frac{\d}{\d x}, \d_1 = \frac{\d}{\d t}\}$  and $\Pi = \{\d_1 \}$. Let  
$k_0$ be a differentially closed $\Pi$-field containing the coefficients of powers of $x$
appearing in the entries of $A$ and let $k = k_0(x)$ be the $\Delta$-field gotten by extending
$\d_1$ via $\d_1(x) = 0$ and defining $\d_0$ to be zero on $k_0$ and $\d_1(x) = 1$. 
\begin{cor} Let $k_0, k, A(t,x)$ be as above and let $K$ be the  PV-extension associated with
\begin{eqnarray}
\frac{\d Y}{\d x} = A(x,t) Y \ . \label{eqn5}
\end{eqnarray}
Then, either
\begin{enumerate}
\item $\PGal(K/k) = \SL_2(k_0)$, or
\item  $\PGal(K/k)$ contains a solvable subgroup of finite index and $K$ is a parameterized
liouvillian extension of $k$, or
\item  equation (\ref{eqn5}) is isomonodromic on $D'\times \calU$ where $D'$ is a punctured subset
of $D$ and $U$ is an open subset of $D(\tau_0,\rho_0)$.
 \end{enumerate}
\end{cor}

We can also state a result similar to  Proposition~\ref{prop2x2} for parameterized linear
equations having more than one parameter. 
We recall that if $k_0$ is a $\Pi =
\{\d_1, \ldots \d_m\}$-field, we denote by $\DX$ the $k_0$-vector space of derivations spanned by
$\Pi$.
\begin{prop} Let $k$ be a $\Delta=\{\d_0, \ldots \d_m\}$-field, assume that $k_0 =C^0_k$ is a
differentially closed $\Pi=\{\d_1, \ldots \d_m\}$-field and let $C = C^\Delta_k$. Let
$A\in\sl_2(k)$ and let $K$ be the PPV-extension corresponding to the differential equation 
\begin{eqnarray*}
\d_0Y = AY \ . 
\end{eqnarray*}
Then, either
\begin{enumerate}
\item $\PGal(K/k) = \SL_2(k_0)$, or
\item  $\PGal(K/k)$ contains a solvable subgroup of finite index and $K$ is a parameterized
liouvillian extension of $k$, or
\item  $\PGal(K/k)$is  a proper Zariski-dense subgroup of $\SL_2(k_0)$ and there exist 
\begin{enumerate}
  
\item   a commuting  $k_0$-basis  $ \{\d_1^\prime, \ldots ,\d_m^\prime \}\mbox{ of } \DX$ and  
\item an integer $r, 1\leq r\leq m$ and elements $B_i \in \gl_2(k), i = 1, \ldots r$,
\end{enumerate} such that the system
\begin{eqnarray*}
\d_0 Y &= & AY\\
\d_1^\prime Y & = & B_1Y\\
  & \vdots & \\
 \d_r^\prime Y & = & B_r Y
 \end{eqnarray*}
 is an integrable system.
 \end{enumerate}
\end{prop} 
\begin{proof} The proof begins in the same way as that for Proposition~\ref{prop2x2}. and Cases 1
and 2 are the same.  If neither of these hold, then $\PGal(K/k)$ is a proper Zariski-dense subgroup
of $\SL_2(k_0)$ and so by Proposition~\ref{sln}, there exist commuting derivations $\Delta' =
\{\d_1^\prime,
\ldots ,\d_r^\prime\}
\subset
\DX$ such that $\PGal(K/k)$ is conjugate to $\SL_2(C_0^{\Delta'})$.  We may assume that the
$\d_i^\prime$ are $k_0$ independent. Proposition 7 of Chapter 0 of \cite{kolchin_groups} implies
that we can extend $\Delta'$ to a commuting basis of $\DX$. After conjugation by an element $B
\in
\GL_2(k)$, we can assume that the PPV-group is  $\SL_2(C_0^{\Delta'})$. A calculation shows that
$(\d_i^\prime Y)Y^{-1}$ is left invariant by this group for $i = 1, \ldots ,r$ and the conclusion
follows.
\end{proof}

The third case of the previous proposition can be stated informally as: {\em After a change of
variables in the parameter space, the parameterized differential equation is completely integrable
with respect to $x$ and a subset of the new parameters}.

\section{Inverse Problem}\label{inversesec} The general inverse problem can be stated as: {\em Given a differential
field, which groups can occur as Galois groups of PPV-extensions of this field?}  We have no
definitive results but will give two examples in this section. 
\begin{example}\label{inverseex1}{\em Let $k$ be a $\Delta = \{\d_0, \d_1\}$-field with $k_0 = C^{\{\d_0\}}_k$
differentially closed and $k = k_0(x), \d_0(x) = 1,\d_1(x) = 0$.  We wish to know: {\it Which
subgroups $G$ of $\Ga(k_0)$ are Galois groups of PPV-extensions of $k$?} The answer is that all proper differential algebraic subgroups 
of $\Ga(k_0)$ appear in this way but $\Ga(k_0)$ itself cannot be the  Galois group of a PPV-extension $K$  of $k$.\\[0.1in]
We begin by showing that $\Ga(k_0)$ cannot be the Galois group of a PPV-extension $K$  of $k$. In Section~\ref{ppvringsandtorsors}, we show that $K$ is the differential 
function field of a $G$-torsor.
If $G = \Ga(k_0)$, then the Corollary to Theorem 4 of Chapter VII.3 of \cite{kolchin_groups} implies that this torsor is trivial and so
$K = k\langle z \rangle_\Delta$ where $\sigma(z) = z + c_\sigma$ for all $\sigma \in \Ga(k_0)$.  This further implies that
$\d_0(z) =a$ for some $a\in k$. Since $k = k_0(x)$ and $k_0$ is algebraically closed, we
may write 
\[
a = P(x) + \sum_{i=1}^r (\sum_{j=1}^{s_i} \frac{b_{i,j}}{(x-c_i)^j})
\]
where $P(x)$ is a polynomial with coefficients in $k_0$ and the $b_{i,j}, c_i \in k_0$. Furthermore,
there exists an element $R(x) \in k$ such that
\[
\d_0(t-R(x)) = \sum_{i=1}^r  \frac{b_{i,1}}{(x-c_i)}
\]
so after such a change, we may assume that 
\[
a = \sum_{i=1}^r \frac{b_i}{x-c_i}
\]
for some $b_i,c_i \in k_0$. \\
[0.1in]
 We shall show that the Galois group of $K$ over $k$ is  
\[
\{ z \in k_0 \ | \ L(z) = 0\}
\]
where $L$ is the linear differential equation in $k[\d_1]$ whose solution space is spanned (over
$C$) by
$\{b_{1}, b_{2}, \dots , b_{r}\}$.  In particular, {\em the group $\Ga(k_0)$ is not a Galois group
of a PPV-extension of
$k$}.\\
[0.1in] To do this
 form a new PPV-extension $F = k\langle z_1,
\ldots , z_r\rangle_\Delta$ where 
$\d_0 z_ i = \frac{1}{x-c_i}$. Clearly, there exists an element $w =\sum_{i=1}^r b_iz_i\in F$ such
that $\d_0 w = a$. Therefore we can can consider $K$ as a subfield of $F$.   A calculation shows
that $\d_0(\d_1 z_i - \frac{\d_1c_i}{x-c_i}) = 0$ so $\d_1z_i \in k$.  Therefore
Proposition~\ref{propint} implies that the PPV-group $\PGal(F/k)$ is of the form $G(C)$ for some
linear algebraic group $G$ and that $F$ is a PV-extension of $k$. The Kolchin-Ostrowski Theorem
(\cite{DAAG}, p.407) implies that the elements $z_i$ are algebraically independent over $k$. The
PPV-group $\PGal(F/k)$ is clearly a subgroup of $\Ga(C)^r$ and since the transcendence degree of $F$
over $k$ must equal the dimension of this group, we have $\PGal(F/k)=\Ga(C)^r$.\\
[0.1in]
For $\sigma= (d_1, \ldots , d_r) \in \Ga(C)^r = \PGal(F/k), \ \sigma(w) = w + \sum_{i=1}^r d_i b_i$.
The Galois theory implies that restricting elements of $\PGal(F/k)$ to $K$ yields a surjective
homomorphism onto $\PGal(K/k)$, so we can identify $\PGal(K/k)$ with the $C$-span of the $b_i$. 
Therefore $\PGal(K/k)$ has the desired form.\\[0.1in]
We now show that any proper differential algebraic subgroup $H$ of $\Ga(k_0)$ is the PPV-group of a PPV-extension of $k$. As stated in Section~\ref{multgps}.
$H =\{a \in \Ga(k_0) \ | \ L(a) = 0\}$ for some linear differential operator $L$ with coefficients in $k_0$.  Since $k_0$ is differentially closed,
there exist $b_1, \ldots , b_m \in k_0$ linearly independent over $C= C_k^{\delta}$ that span the solution space of $L(Y) = 0$.
Let 
\[ a = \sum_{i=1}^m \frac{b_i}{x-i} \ . \] 
The calculation above shows that the PPV-group of the PPV-extension of $k$ for $\d_0 y = a$ is $H$.\hfill \QED }\end{example}
The previous example leads to the question: Find a $\Delta$-field $k$ such that $\Ga(k_0)$ is a
Galois group of a PPV-extension of $k$. We do this in the next example.
\begin{example}\label{gammafnc}{\em  Let $k$ be a $\Delta = \{\d_0, \d_1\}$-field with $k_0 = C^{\{\d_0\}}_k$
differentially closed and $k = k_0(x, \log x, x^{t-1}e^{-x}), \d_0(x) = 1,\d_0(\log x) =\frac{1}{x},
\d_0(x^{t-1}e^{-x}) = \frac{t-x-1}{x}x^{t-1}e^{-x}, \d_1(x) = \d_1(\log x) = 0,
\d_1(x^{t-1}e^{-x}) = (\log x) x^{t-1}e^{-x} $. Consider the differential equation 
\[
\d_0 y = x^{t-1}e^{-x}
\]
and let $K$ be the PPV-extension of $k$ for this equation. We may write $K
 =k\langle \gamma\rangle_\Delta$, where $\gamma$ satisfies the above equation ($\gamma$ is known as
the {\em incomplete Gamma function}).  We have that $K = k(
\gamma, \d_1\gamma, \d_1^2\gamma, \ldots)$.   In \cite{JRR}, the authors show that $\gamma,
\d_1\gamma, \d_1^2\gamma,
\ldots$ are algebraically independent over $k$.  Therefore,  for any $c \in \Ga(k_0), \ 
\d_1^i \gamma \mapsto \d_1^i\gamma +\d_1^ic, \ i=0, 1, \ldots $ defines an element
of $\Gal(K/k)$. Therefore  $\Gal(K/k) = \Ga(k_0)$.\\
[0.1in]
Over $ k_0(x)$, $\gamma$ satisfies
\[
\frac{\d^2\gamma}{\d x^2}- \frac{t-1-x}{x} \  \frac{\d\gamma}{\d x} = 0
\]
and one can furthermore show that the PPV-group over $k_0(x)$ of this latter equation is 
\begin{eqnarray*}
H& =& \{ \left( 
\begin{array}{cc} 1&a\\
0&b 
\end{array}\right) \ | \ a \in k_0, b\in 
k_0^*,\d_1(\frac{\d_tb}{b} )=0\}\\
  &=& \Ga(k_0) \rtimes\Gm^{\d_1},  
  \end{eqnarray*} 
  where $\Gm^{\d_1} = \{ b\in k_0^* \ | \ \d_1(\frac{\d_1 b}{b}) = 0\}$.\hfill \QED }
\end{example}

\begin{remark}\label{gamma}{\em  We can use the previous example to exhibit two equations $\d_0 Y = A_1Y$ and $\d_0Y = A_2Y$ having the same
PPV-extension $K$ of $k$ but such that $K_{A_1}^{\rm PV} \neq K_{A_2}^{\rm PV}$ and that these 
latter PV-extensions have different PV-groups ({\em cf.}, Remark~\ref{KPVremark}).
Let $k$ and $\gamma$  be as in the above example and let 
\[A_1 = \left(\begin{array}{cc} 0 & x^{t-1}e^{-x} \\
0 & 0 \end{array}\right) \hspace{1in} A_2 = \left(\begin{array}{ccc} 0 & x^{t-1}e^{-x} & (\log x) x^{t-1}e^{-x}\\
0 & 0& 0  \\ 0 & 0 & 0 \end{array}\right) \ . \]
We have that 
\[ Z_1 = \left(\begin{array}{cc} 1 & \gamma \\
0 & 1 \end{array}\right) \hspace{1in} Z_2 = \left(\begin{array}{ccc} 1 & \gamma & \d_1(\gamma)\\
0 & 1& 0  \\ 0 & 0 & 1 \end{array}\right) \]
satisfy $\d_0 Z_1 = A_1 Z_1$ and $\d_0 Z_2 = A_2 Z_2$. $K$ is the PPV-extension associated with either equation and the 
Galois group $\PGal(K/k)$ is $\Ga(k_0)$. We have that $K^{\rm PV}_{A_1} = k(\gamma) \neq K^{\rm PV}_{A_2} = k(\gamma, \d_0\gamma)$
since $\gamma$ and $\d_1 \gamma $ are algebraically independent over $k$. With respect to the first equation, $\PGal(K/k)$ is represented in $\GL_2(k_0)$ as
\[ \{\left(\begin{array}{cc} 1 & c \\
0 & 1 \end{array}\right) \ | \ c \in k_0\}\]
and with respect to the second equation $\PGal(K/k)$ is represented in $\GL_3(k_0)$ as
\[ \{ \left(\begin{array}{ccc} 1 & c & \d_1(c)\\
0 & 1& 0  \\ 0 & 0 & 1 \end{array}\right) \ | \ c \in k_0 \} \ . \]
The image of $\Ga(k_0)$ in $\GL_2(k_0)$ is Zariski-closed while the Zariski closure of the image of $\Ga(k_0)$ in $\GL_3(k_0)$ is
\[ \{ \left(\begin{array}{ccc} 1 & c & d\\
0 & 1& 0  \\ 0 & 0 & 1 \end{array}\right) \ | \ c, d \in k_0 \} \ . \]
As algebraic groups, the first group is just $\Ga(k_0)$ and the second is $\Ga(k_0) \times \Ga(k_0)$.
}\end{remark}
\section{Final Comments}\label{finalsec}
\noindent{\bf Other Galois Theories}  In \cite{pillaygalois3}, Pillay proposes a Galois theory that
extends Kolchin's Galois theory of strongly normal extensions.  We will explain the connections 
to our results.\\[0.1in]
Let $k$ be a differential field and $K$ a Picard-Vessiot extension of $k$.  $K$ has the following property: for any
differential extension $E$ of $K$ and any differential $k$-isomorphism $\phi$ 
of $K$ into $E$, we have that $\phi(K)\cdot C =  K\cdot C$, where $C$ is the field of constants
of $E$.  Kolchin has shown (\cite{DAAG}, Chapter VI) this is the key property for developing a Galois theory.
In particular, he defines a finitely generated differential field extension $K$ of $k$ to be {\em strongly normal} if for any
differential extension $E$ of $K$ and any differential $k$-isomorphism of $K$ into $E$ we have 
that  \begin{enumerate}
\item[(1)] $\phi(K)\langle C\rangle = K\langle C\rangle$,  where $C$ are the constants of $E$ and
\item [(2)]$\phi$  leaves each of the constants of $K$ fixed.
\end{enumerate} 
For such fields, Kolchin shows that the differential Galois group of $K$ over $k$ has the structure of an 
algebraic group and that the usual Galois correspondence holds. \\[0.1in]
In \cite{pillaygalois1,pillaygalois2,pillaygalois3, pillaygalois4} Pillay considers {\em ordinary} differential fields and
generalizes this theory.  The key observation is that the condition
(1) can be restated as
\begin{enumerate}
\item[(1')]$\phi(K)\langle X(E)\rangle = K\langle X(E)\rangle$, 
\end{enumerate} 
where $X$ is the differential algebraic variety defined by the equation $\d Y = 0$ and $X(E)$ are the $E$-points of $X$.
For $X$, any differential algebraic variety defined over $k$ (or more generally, any Kolchin-constructible set), Pillay defines
a differential extension $K$ to be an {\em $X$-strongly normal} extension of $k$ if for any
differential extension $E$ of $K$ and any differential $k$-isomorphism of $K$ into $E$ we have that equation (1') holds
and that (2) is replaced by  technical (but important) other  conditions. Pillay then uses model theoretic tools to show that for these extensions, the Galois group is a 
finite dimensional differential algebraic group (note that in the PPV-theory, infinite dimensional differential algebraic
groups can occur, {\em e.g.}, $\Ga$).  The finite dimensionality results from the fact that the underlying differential fields are ordinary
differential fields and that finite sets of elements in the differential closure of an ordinary differential field generate fields of finite transcendence degree (a fact that is no longer true 
for partial differential fields).  Because of this, Pillay was able to recast his theory in \cite{pillaygalois4} in the language of subvarieties of certain jet
spaces.  If one generalizes Pillay's definition of strongly normal to allow {\em partial} differential fields with derivations $\Delta$ and takes for $X$ the 
differential algebraic variety defined by $\{ \d_Y = 0 \ | \ \d \in \Pi\}$ where $\Pi \subset \Delta$, then this definition
would include PPV-extensions.  Presumably the techniques  of \cite{pillaygalois3} can be used to prove many of these results as well.  Nonetheless,
we feel that a description of the complete situation for PPV-fields is sufficiently self contained as to warrant an independent exposition. \\[0.1in]
Landesman \cite{Landesman} has been generalizing Kolchin's Galois 
theory of strongly normal extensions  to differential fields having a designated subset of derivations acting as parametric 
derivations. When this is complete, many of  our results should follow as a special case of his results.\\[0.1in]
Umemura~\cite{u1,u2,u3,u4,u5,u6} has proposed a Galois theory for general nonlinear differential equations. Instead of
Galois groups, he uses Lie algebras to measure the symmetries of differential fields.  Malgrange~\cite{malgrange_galois1, malgrange_galois2}
has proposed a Galois theory of differential equations where the role of the Galois group is taken by certain groupoids.  Both
Umemura and Malgrange have indicated to us that their theories can analyze parameterized differential equations as well. \\[0.2in]
\noindent{\bf Future Directions} There are many questions suggested by the results presented here and we indicate a few of these.
\begin{enumerate}
\item Deligne \cite{deligne_milne, deligne_tannakian} (see also \cite{PuSi2003}) has shown that the usual 
Picard-Vessiot theory can be presented in the language of Tannakian Categories. Can one characterize in a similar way the 
category of representations  
of linear differential algebraic groups and use this to develop the Parameterized Picard-Vessiot theory?
\item How does the parameterized monodromy sit inside the parameterized Picard-Vessiot groups? To what extent can one
extend Ramis' characterization of the local Galois groups to the parameterized case?
\item Can one develop algorithms to determine the Galois groups of parameterized linear differential equations. 
Sit \cite{sit75} has classified the differential algebraic subgroups of $\SL_2$. Can this classification be used to 
calculate Galois groups of second order paramterized differential equations in analogy to 
Kovacic's algorithm for second order linear differential equations?
\item Characterize those linear differential algebraic groups that appear as Galois groups of $k_0(x)$ where $k_0$ is as in Example~\ref{inverseex1}.
\end{enumerate}

\section{Appendix}
In this Appendix, we present proofs of results that imply Theorem~\ref{PPVthm} and Theorem~\ref{liouvillian}.
In Section 3, Theorem~\ref{PPVthm} is stated for a parameterized system of {\em ordinary} linear differential equations but it is no harder to
prove an analogous result for parameterized integrable systems of linear partial differential equations and we do this in this appendix.
The
first section contains a discussion of {\em constrained extensions}, a concept needed in the proof
of the existence of PPV-extensions. In the next three sections, we prove results that simultaneously imply Theorem~\ref{PVthm} and
Theorem~\ref{PPVthm}. The proofs are almost, word-for-word, the same as the proofs of the
corresponding result for PV-extensions (\cite{PuSi2002}, Ch.~1) once one has taken into account the need for subfields of constants 
to be differentially closed.  Nonetheless  we include the proofs for the
convenience of the reader. The final section contains a proof of Theorem~\ref{liouvillian}.
\subsection{Constrained Extensions} Before turning to the proof of Theorem~\ref{PPVthm}, we shall
need some more facts concerning 
differentially closed fields (see Definition~\ref{diffclosed}). If $k\subset K$ are $\Delta$-fields and $\eta = 
(\eta_1, \ldots, \eta_r) \in K^r$, we denote
by $k\{\eta\}_\Delta$ (resp.~$k\langle \eta \rangle_\Delta$) the $\Delta$-ring (resp.~$\Delta$-field) generated by $k$ and $\eta_1, \ldots , \eta_r$,
that is, the ring (resp.~field) generated by $k$ and all the derivatives of the $\eta_i$.  We shall denote by $k\{y_1, \ldots , y_n\}_\Delta$ the ring
of differential polynomials in $n$ variables over $k$ ({\em cf.,} Section 3). A $k$-$\Delta$-isomorphism of $k\{\eta\}_\Delta$ is a $k$-isomorphism $\sigma$ such that
$\sigma \d = \d \sigma$ for all $\d \in \Delta$.
\begin{defin}\label{constraineddef}(\cite{DAAG}, Ch.~III.10; \cite{constrained}) Let $k\subset K$ be
$\Delta$-fields.  
\begin{enumerate}
\item We say that a finite family of elements $\eta =(\eta_1, \ldots ,\eta_r) \subset K^r$ is {\em
constrained} over $k$ if there exist
 differential polynomials $P_1, \ldots , P_s, Q  \in k\{y_1, \ldots ,y_r\}_\Delta$ such that 
\begin{enumerate}
\item $P_1(\eta_1, \ldots , \eta_r) = \ldots = P_1(\eta_1, \ldots , \eta_r) = 0$ and $Q(\eta_1,
\ldots , \eta_r)\neq 0$, and 
\item for any $\Delta$-field $E, k \subset E$, if $(\zeta_1, \ldots , \zeta_r) \in E^r$ and $P_1(\zeta_1,
\ldots , \zeta_r) = \ldots = P_1(\zeta_1, \ldots , \zeta_r) = 0$ 
and $Q(\zeta_1, \ldots , \zeta_r)\neq 0$, then the map $\eta_i \mapsto \zeta_i$ induces a
$k$-$\Delta$-isomorphism of $k\{\eta_1, \ldots , \eta_r\}_\Delta$ with $k\{\zeta_1, \ldots
,\zeta_r\}_\Delta$.
\end{enumerate}
We say that $Q$ is the {\em constraint} of $\eta$ over $k$.
\item We say $K$ is a {\em constrained} extension of $k$ if every finite family of elements of $K$
is constrained over $k$.
\item We say $k$ is {\em constrainedly closed} if $k$ has no proper constrained extensions. 
\end{enumerate}
\end{defin} 
The following Proposition contains the facts that we will use:
\begin{prop}\label{constrainedfacts} Let $k \subset K$ be $\Delta$-fields and $\eta \in K^r$
\begin{enumerate}
\item $\eta$ is {\em constrained} over $k$ with constraint $Q$ if and only
if $k\{\eta, 1/Q(\eta)\}_\Delta$ is a simple $\Delta$-ring, {\em i.e.} a $\Delta$-ring with no
nontrivial $\Delta$-ideals. 
\item If $\eta$ is constrained over $k$ and $K=k\langle \eta\rangle_\Delta$, 
then any finite set of
elements of $K$ is constrained over $k$, that is, $K$ is a constrained extension of $k$.
\item $K$ is differentially closed if and only if it is constrainedly closed.
\item Every differential field has a constrainedly closed extension.
\end{enumerate}
\end{prop} One can find the proofs of these in \cite{constrained}, where Kolchin uses the term
constrainedly closed instead of differentially closed. Proofs also 
can be found in \cite{mcgrail} where the author uses a model theoretic approach. 
Item 1.~follows from
the fact that any maximal $\Delta$-ideal in  a ring containing  $\QX$ is prime (\cite{DAAG}, Ch.~I.2,
Exercise 3 or \cite{PuSi2003}, Lemma~1.17.1) and that for any radical differential ideal $I$ in
$k\{y_1, \ldots , y_r\}_\Delta$ there exist differential polynomials $P_1, \ldots ,P_s$ such that
$I$ is the smallest radical differential ideal containing $P_1, \ldots ,P_s$ (the Ritt-Raudenbusch
Theorem \cite{DAAG}, Ch.~III.4). Item 2.~is fairly deep and is essentially equivalent to the fact
that the projection of a Kolchin-constructible set (an element in the boolean algebra generated by
Kolchin-closed sets)
 is Kolchin-constructible. Items 3.~and 4.~require some effort but are not too difficult to
prove. Generalizations to fields with noncommuting derivations can be found in \cite{yaffe} and
\cite{pierce}. \\[0.1in]
In the usual Picard-Vessiot theory, one needs the following key fact: Let $k$ be a differential field with algebraically closed
subfield of constants $C$. If $R$ is a simple differential ring, finitely generated over $k$, 
then any constant of $R$ is in $C$ (Lemma 1.17, \cite{PuSi2003}). The following result
generalizes this fact and plays a similar role in the Parameterized Picard-Vessiot Theory.
Recall that if $k$ is a 
 $\Delta = \{\d_0, \ldots ,\d_m\}$-field
 and $\Lambda \subset \Delta$, we denote by $C_k^{\Lambda}$ the set $\{c \in k \ | \ \d c = 0 \ 
\forall \d \in \Lambda\}$. One sees that $C_k^\Lambda$ is a $\Pi = \Delta\backslash
\Lambda$-field. 
 \begin{lem}\label{constants} Let $k \subset K$ be $\Delta$-fields,  $\Lambda \subset \Delta, $ and $ \Pi = \Delta
\backslash \Lambda$. Assume that $C_k^\Lambda$ is $\Pi$-differentially closed. 
   If  $K$ is a $\Delta$-constrained
extension of $k$, then $C^{\Lambda}_K=C^{\Lambda}_k$.
 \end{lem}
 \begin{proof} Let $\eta \in C^{\Lambda}_K$.  Since $K$ is a $\Delta$-constrained extension of $k$, there
exist $P_1, \ldots , P_s, Q  \in k\{y\}_\Delta$ satisfying the conditions of
Definition~\ref{constraineddef} with respect to $\eta$ and $k$.  We will first show that there exist
$P_1, \ldots , P_s, Q  \in C_k^{\Lambda}\{y\}_\Delta$ satisfying the conditions of
Definition~\ref{constraineddef} with respect to $\eta$ and $k$.\\[0.1in]
 Let $\{\beta_i\}_{i \in I}$ be a
$C^{\Lambda}_k$-basis of $k$. 
Let $R \in k\{y\}_\Delta$ and write $R = \sum R_{i} \beta_i
$ where each $R_i
\in C^{\Lambda}_k\{y\}_\Delta$. Since linear independence over constants is
preserved when one goes to  extension fields (\cite{DAAG}, Ch.~II.1),   for any
differential $\Delta$-extension $E$ of $k$ and $\zeta \in C_E^{\Lambda}$, we have that 
$R(\zeta) = 0$ if and only if all $R_i(\zeta) = 0$ for all $i$. If we write $P_j = \sum
P_{i,j}\beta_i, Q =
\sum Q_i\beta_i$ then there is some $i_0$ such that $\eta$ satisfies $\{P_{i,j} = 0\}, Q_{i_0} \neq
0$ and that for any $\zeta \in C_E^{\Lambda}$ that satisfies this system, the map $\eta \mapsto \zeta$
induces a $\Delta$ isomorphism of $k\langle \eta \rangle_\Delta$ and 
$k\langle \zeta \rangle_\Delta$.\\[0.1in]
We therefore may assume that there
exist $P_1, \ldots , P_s, Q  \in C^{\Lambda}_k\{y\}_\Delta$ satisfying the conditions of
Definition~\ref{constraineddef} with respect to $\eta$ and $k$. We now show that there exist 
$\tilde{P}_1, \ldots , \tilde{P}_s, \tilde{Q}$ 
in the smaller differential polynomial ring $C^{\Lambda}_k\{y\}_{\Pi}$ satisfying:
If $E$ is a $\Delta$-extension of $k$ and $\zeta \in C_E^{\Lambda}$ satisfies $\tilde{P}_1(\zeta) =\ldots \tilde{P}_s(\zeta)=0, 
\tilde{Q}(\zeta) \neq 0$ then there is a $k$-$\Delta$-isomorphism of $k\langle \eta\rangle_\Delta$ and $k\langle \zeta\rangle_\Delta$
mapping $\eta \mapsto \zeta$.  To do this, note that any $P \in 
C^{\Lambda}_k\{y\}_\Delta$ is a $C^{\Lambda}_k$-linear combination of monomials
that are
products of terms of the form $\d_0^{i_0} \cdots \d_m^{i_m}y$.  We
denote by $\tilde{P}$ the differential polynomial resulting from $P$ be deleting any monomial that
contains a term $\d_0^{i_0} \cdots \d_m^{i_m}y_j$ {\em with $i_t >0$ for some $\d_{i_t}\in
\Lambda$}.  Note that for any $\Delta$-extension $E$ of $k$ and $\zeta \in C_E^{\Lambda}$ we
have $P(\zeta) = 0$ if and only if $\tilde{P}(\zeta) = 0$.  Therefore, for any  $\zeta \in
C_E^{\Lambda}$,  if $\tilde{P}_1(\zeta) = \ldots = \tilde{P}_1(\zeta) = 0$ 
and $\tilde{Q}(\zeta)\neq 0$, then the map $\eta \mapsto \zeta$ induces a
$\Delta$-$k$-isomorphism of $k\{\eta\}_\Delta$ with $k\{\zeta\}_\Delta$. \\[0.1in]
We now use the fact that $C_k^\Lambda$ is a $\Pi$-differentially closed field to show that any $\eta \in C_K^\Lambda$ must
already be in $C_k^\Lambda$.  Let $\tilde{P}_1, \ldots , \tilde{P}_s, \tilde{Q} \in 
C^{\Lambda}_k\{y\}_{\Pi}$ be as above. Since $C_k^\Lambda$ is a $\Pi$-differentially closed field and
$\tilde{P}_1= \ldots =\tilde{P}_s=0, \tilde{Q}\neq 0$ has a solution in {\em some} $\Pi$-extension of $C_k^\Lambda$ ({\em e.g.,}
$\eta \in C_K^\Lambda$), this system has a solution $\zeta \in C_k^\Lambda \subset k$. We therefore can conclude that
the map $\eta \mapsto \zeta$ induces a  $\Pi$-$k$-isomorphism from $k\langle \eta\rangle$ to
to $k\langle \zeta \rangle$. Since $\zeta \in k$, we have that $\eta \in k$ and so $\eta \in C_k^\Lambda$.\end{proof}
We note that if $\Pi$ is empty, then $\Pi$-differentially closed is the same as algebraically closed.  In this case the above result yields the 
important fact crucial to the Picard-Vessiot theory mentioned before the lemma.

\subsection{PPV-extensions} In the next three sections, we will develop the theory of PPV-extensions for parameterized integrable systems of linear differential equations. This section is devoted to
showing the existence and uniqueness of these extensions, Section~\ref{subgaloisgroups} 
we show that the Galois group has a natural structure as a linear differential algebraic group and in Section~\ref{ppvringsandtorsors} we show that
a PPV-extension can be associated with a torsor for the Galois group.  As in the usual Picard-Vessiot theory, these results will allow us to give a complete Galois theory (see Theorem~\ref{completegalois}).\\[0.1in]
In this and the next three sections, we will make  the following conventions. We let $k$ be a $\Delta$-differential field. 
We designate a nonempty subset $\Lambda =\{\ \d_0, \ldots , \d_r\} \subset \Delta$ and 
consider a system of linear differential equations 
\begin{eqnarray}
\d_0 Y&=& A_0Y \label{intsys2}\\
\d_1 Y  &=&  A_0Y\notag\\
   &\vdots&  \notag\\
  \d_r Y  &= & A_r Y \notag
  \end{eqnarray}
  where the $A_i \in \M_n(k)$, the set of $n\times n$ matrices with entries in $k$, such that
 \begin{eqnarray}
\d_iA_j - \d_jA_i & = & [A_i, A_j] 
\end{eqnarray}
We denote by $\Pi$ the set $\Delta\backslash\Lambda$. One sees that the derivations of $\Pi$ leave the field $C_k^\Delta$ invariant
and we shall think of this latter field as a $\Pi$-field.  Throughout the next sections, we shall assume that {\em 
$C = C_k^\Delta$ is a $\Pi$-differentially closed differential field}. The set $\Lambda$ corresponds to derivations used in the \underline{l}inear differential equations and 
$\Pi$ corresponds to the \underline{p}arametric derivations.Throughout the first part of this paper $\Delta$ was $\{\d_0, \ldots ,\d_m\},
\Lambda = \{\d_0\}, \mbox{ and } \Pi = \{\d_1, \ldots, \d_m\}$. We now turn to a definition.

\begin{defin} \begin{enumerate}
\item A {\em Parameterized Picard-Vessiot ring} (PPV-ring) over  $k$ for
the equations (\ref{intsys2}) is a $\Delta$-ring $R$ containing $k$ satisfying:
\begin{enumerate}
\item $R$ is a $\Delta$-simple $\Delta$-ring.
\item There exists a matrix $ Z \in \GL_n(R)$ such that $\d_i Z = A_iZ$ for all $\d_i \in \Lambda$.
\item $R$ is generated, as a $\Delta$-ring over $k$, by the entries of $Z$ and $1/\det(Z)$, {\em i.e.}, 
$R=k\{Z,1/\det(Z)\}_\Delta$.
\end{enumerate}
\item A 
{\em Parameterized Picard-Vessiot extension of $k$  (PPV-extension of $k$)}  for the equations
(\ref{intsys2}) is a $\Delta$-field $K$ satisfying 
\begin{enumerate}
\item $k \subset K$.
\item There exists a matrix $ Z \in \GL_n(K)$ such that $\d_i Z = A_iZ$ for all $\d_i \in \Lambda$.
\item $C^\Lambda_K = C^\Lambda_k$, {\em i.e.,} the $\Lambda$-constants of $K$ coincide with the $\Lambda$-constants of
$k$.
\end{enumerate}
\item
The group $\PGal(K/k) = \{\sigma : K\rightarrow K \ | \ \sigma \mbox{ is a }
k\mbox{-automorphism such that } \sigma \d = \d \sigma \ \forall  \d \in \Delta \}$ is called
the {\em Parameterized Picard-Vessiot Group (PPV-Group)} associated with the PPV-extension $K$ of
$k$.\end{enumerate}
\end{defin}
Note that when $\Delta = \Lambda, \Pi = \emptyset$ these definitions give us the corresponding definitions in the usual Picard-Vessiot theory.\\[0.1in]
Our goal in the next three sections is to prove results that will yield the following generalization of both 
Theorem~\ref{PVthm} (when $\Delta = \Lambda) 
$ and
Theorem~\ref{PPVthm} (when $\Delta = \{\d_0, \d_1, \ldots , \d_m\}$ and  $\Lambda = \{\d_0\})$.
\begin{thm}\label{completegalois} 
\begin{enumerate}
\item There exists a PPV-extension $K$ of $k$ associated with (\ref{intsys2}) and this is unique up
to $\Delta$-$k$-isomorphism.
\item The PPV-Group $\PGal(K/k)$ equals $G(C_k^\Lambda)$, where $G$ is a linear $\Pi$-differential algebraic group defined over $C_k^\Lambda$.
\item The  map that sends any $\Delta$-subfield $F, k \subset F \subset K$, to the group
$\PGal(K/F)$ is a bijection between  $\Delta$-subfields of $K$ containing $k$ and $\Pi$-Kolchin closed
subgroups of $\PGal(K/k)$. Its inverse is given by the map that sends a $\Pi$-Kolchin closed group $H$ to
the field $\{z \in K \ | \ \sigma(z) = z \mbox{ for all } \sigma \in H\}$.
\item A $\Pi$-Kolchin closed subgroup $H$ of $\PGal(K/k)$ is a normal subgroup of $\PGal(K/k)$ if and
only if the field $K^H$ is left set-wise invariant by   $\PGal(K/k)$. If this is the case, the 
map $\PGal(K/k) \rightarrow  \PGal(K^H/k)$ is surjective with kernel $H$ and
$K^H$ is a PPV-extension of $k$ with PPV-group isomorphic to $\PGal(K/k)/H$. Conversely, if $F$ is
a differential subfield of $K$ containing $k$ and $F$ is a PPV-extension of $k$, then $\PGal(F/K)$
is a normal $\Pi$-Kolchin closed subgroup of $\PGal(K/k)$. 
\end{enumerate}
\end{thm}
  We shall show in this section that
PPV-extensions for (\ref{intsys2}) exist and are unique up to $\Delta$-$k$-isomorphism and that every PPV-extension $K$ of $k$ is 
the quotient field of a PPV-ring (and therefore is also unique up to $\Delta$-$k$-isomorphism.  We begin with
\begin{prop}\label{propPPVquotient} 
\begin{enumerate}
\item There exists a PPV-ring $R$ for (\ref{intsys}) and it is an integral domain.
\item The field of $\Lambda$-constants $C^\Lambda_K$ of the quotient field $K$ of a PPV-ring over $k$ is $C^\Lambda_k$.
\item Any two PPV-rings for this system are $k$-isomorphic as $\Delta$-rings.  
\end{enumerate}
\end{prop}

\begin{proof} 1.  Let $(Y_{i,j})$ denote an $n \times n$ matrix of $\Pi$-indeterminates and let
``det'' denote the determinant of
$(Y_{i,j})$.  We denote by $k\{Y_{1,1}, \ldots, Y_{n.n}, 1/\det \}_{\Pi}$ the $\Pi$-differential polynomial ring in the variables
$\{Y_{i,j}\}$ localized at $\det$. We can make this ring into a $\Delta$-ring by setting 
$(\d_k Y_{i,j}) = A_k (Y_{i,j})$ for all $\d_k \in \Lambda$ and using the fact that $\d_k\d_l
= \d_l \d_k$ for all $\d_k, \d_l \in \Delta$.  Let $p$ be a maximal $\Delta$-ideal in $R$. One then sees that $R/{p}$ is a PPV-ring for
the equation. Since maximal differential ideals are prime,$R$ is an integral domain.\\[0.1in]
2.  Let $R = k \{ Z, 1/\det(Z)\}_\Delta$. Since this is a simple differential ring,
Proposition~\ref{constrainedfacts}.1 implies 
 that $Z$ is constrained over $k$ with constraint $\det$.  Statement 
 2.~of Proposition~\ref{constrainedfacts} implies that the quotient field of $R$ is a
$\Delta$-constrained extension of $k$.  Lemma~\ref{constants} implies that $C_K^{\Lambda} = V_k^{\Lambda}$.\\[0.1in]
3. Let $R_1, R_2$ denote two PPV-rings for the system.  Let $Z_1, Z_2$ be the two fundamental
matrices.  Consider the $\Delta$-ring $R_1\otimes_{k} R_2$ with derivations $\d_i(r_1\otimes
r_2) = \d_ir_1 \otimes r_2 + r_1 \otimes \d_ir_2$.  Let ${p}$ 
be a maximal $\Delta$-ideal in $R_1\otimes_k R_2$ and let $R_3 = R_1\otimes_{k} R_2/{p}$. 
The obvious maps $\phi_i:R_i \rightarrow R_1\otimes_{k} R_2$ are $\Delta$-homomorphisms and,
since the $R_i$ are simple, the homomorphisms $\phi_i$ are injective.  The image of each $\phi_i$ is
differentially generated 
by the entries of $\phi_i(Z_i)$ and $\det(\phi(Z_i^{-1}))$.  The matrices $\phi_1(Z_1)$  and
$\phi_2(Z_2)$ are fundamental matrices in $R_3$ of 
the differential equation. Since $R_3$ is simple, the previous result implies that $C_k^\Lambda$ is the ring of $\Lambda$-constants of
$R_3$.  Therefore 
$\phi_1(Z_1) = \phi_2(Z_2)D$ for some $D \in \GL_n(C^{\Lambda}_k)$.  Therefore $\phi_1(R_1) = \phi_2(R_2)$
and so $R_1$ and $R_2$ are isomorphic.
\end{proof}
Conclusion 2.~of the above proposition shows that the field of fractions of a PPV-ring is a
PPV-field.  We now show that a PPV-field for an equation 
is the field of fractions of a PPV-ring for the equation.
\begin{prop}\label{propPPVfield} Let $K$ be a PPV-extension field of $k$ for the system
(\ref{intsys2}), let $Z \in \GL_n(K)$ satisfy $\d_i(Z) = A_iZ$ for all $\d_i \in \Lambda$ and let $\det = \det(Z)$. 
\begin{enumerate}
\item The
$\Delta$-ring $k\{Z, 1/\det\}_\Delta$ is a PPV-extension ring of $k$ for this system. 
\item If $K'$ is another PPV-extension of $k$ for this system then there is a
$k$-$\Delta$-isomorphism of $K$ and $K'$.
\end{enumerate}
\end{prop}
 To simplify notation we shall use
$\frac{1}{{\rm det}}$ to denote the inverse of the determinant of a matrix given
by the context.  For example, $k\{Y_{i,j},\frac{1}{{\rm det}}\}_{\Delta} = 
k\{Y_{i,j},\frac{1}{{\rm det}(Y_{i,j})}\}_{\Delta}$ and $k\{X_{i,j},\frac{1}{{\rm det}}\}_{\Pi} = 
k\{X_{i,j},\frac{1}{{\rm det}(X_{i,j})}\}_{\Pi}$.\\[0.1in]
 As in \cite{PuSi2003}, p.~16, we need a preliminary lemma to prove this proposition. Let $(Y_{i,j})$
be an $n
\times n$ matrix of $\Pi$-differential indeterminates
 and let $\det$ denote the determinant of this matrix. For any $\Pi$-field $k$, we denote by
$k\{Y_{i,j}, 1/\det\}_{\Pi}$
 the $\Pi$-ring of differential polynomials in the $Y_{i,j}$ localized with respect to  $\det$.
If $k$ is, in addition, a $\Delta$-field, the
 derivations $\d \in \Lambda$ can be extended to 
 $k\{Y_{i,j}, 1/\det\}_{\Pi}$ by setting $\d(Y_{i,j}) = 0$ for all $\d \in \Lambda$ and  $i,j$ with $1 \leq i,j
\leq n$.  In this way
  $k\{Y_{i,j}, 1/\det\}_{\Pi}$ may be considered as a $\Delta$-ring. We
 consider $C_k^\Lambda \{Y_{i,j}, 1/\det\}_{\Pi}$ as a $\Pi$-subring of $k\{Y_{i,j},
1/\det\}_{\Pi}$. For any set $I\subset k\{Y_{i,j}, 1/\det\}_{\Pi}$, we denote by $(I)_\Delta$ the $\Delta$-differential ideal
in $k\{Y_{i,j}, 1/\det\}_{\Pi}$ generated by $I$.
 
 \begin{lem}\label{lem2} Using the above notation, the map $I \mapsto (I)_\Delta$ from the set of
 $\Pi$-ideals of $C_k^\Lambda\{Y_{i,j}, 1/\det\}_{\Pi}$ to the set of $\Delta$- ideals of
$k\{Y_{i,j}, 1/\det\}_{\Pi}$ is
 a bijection.  The inverse map is given by $J \mapsto J \cap C_k^\Lambda\{Y_{i,j}, 1/\det\}_{\Pi}$.
 \end{lem}
 \noindent {\bf Proof of Lemma~\ref{lem2}.} If $\calS = \{s_{\alpha}\}_{\alpha \in \calA}$ is a
basis of $k$ over $C_k^\Lambda$, then $\calS$ is a module basis of $k\{Y_{i,j},\frac{1}{{\rm
det}}\}_{\Pi}$ over $C_k^\Lambda\{Y_{i,j},\frac{1}{{\rm det}}\}_{\Pi}$.  Therefore, for any ideal $I$
of $C_k^\Lambda\{Y_{i,j},\frac{1}{{\rm det}}\}_{\Pi}$, one has that $(I)_\Delta \cap
C_k^\Lambda\{Y_{i,j},\frac{1}{{\rm det}}\}_{\Pi} = I$.\\[0.1in] 
We now prove that any ${\Delta}$-differential ideal $J$ of $k\{Y_{i,j},\frac{1}{{\rm
det}}\}_{\Pi}$ is generated by $I := J \cap C_k^\Lambda\{Y_{i,j},\frac{1}{{\rm det}}\}_{\Pi}$.  Let
$\{e_{\beta}\}_{\beta \in \calB}$ be a basis of $C_k^\Lambda\{Y_{i,j},\frac{1}{{\rm det}}\}_{\Pi}$ over
$C_k^\Lambda$.  Any element $f \in J$ can be uniquely written as a finite sum $\sum_{\beta \in \calB}m_{\beta}e_{\beta}$
with the $m_{\beta} \in k$.  By induction on the length, $l(f)$, of $f$ we will show that $f \in
(I)_\Delta$.  When $l(f) = 0,1$, the result is clear.  Assume $l(f) > 1$.  We may suppose that $m_{\beta_1}
= 1$ for some $\beta_1 \in \calB$ and $m_{\beta_2} \in
k\backslash C_k^\Lambda$ for some $\beta_2 \in \calB$.  One then has that, for any $\d \in \Lambda$, $\d f =
\sum_{\beta}\d m_{\beta}e_{\beta}$ has a length smaller than $l(f)$ and so belongs to $(I)_\Delta$. 
Similarly $\d(m_{\beta_2}^{-1}f) \in (I)_\Delta$.   Therefore $\d(m_{\beta_2}^{-1})f =
\d(m_{\beta_2}^{-1}f) - m_{\beta_2}^{-1}\d f \in (I)$. 
Since $C_k^\Lambda$ is the field of $\Lambda$-constants of $k$, one has $ \d_i(m_{\beta_2}^{-1}) \neq 0$ for some $\d_i \in \Lambda$ and so
$f
\in (I)_\Delta $. 
~\hspace{\fill}~$\square$

\vspace{.2in}

\noindent {\bf Proof of Proposition~\ref{propPPVfield}.}
1.~Let $R_0 = k\{X_{i,j}, \frac{1}{det}\}_{\Pi}$ be the 
 ring of $\Pi$-differential polynomials over $k$ and define a $\Delta$-structure on this ring by 
 setting $(\d_iX_{i,j}) = A_i(X_{i,j})$ for all $\d_i \in \Lambda$.  Consider the $\Delta$-rings  
$R_0\subset K\otimes _k R_0=K\{X_{i,j},\frac{1}{{\rm det }}\}_{\Pi}$. Define 
a set of $n^2$ new variables $Y_{i,j}$ by $(X_{i,j})=Z\cdot (Y_{i,j})$. Then 
$K\otimes _k R_0=K\{Y_{i,j},\frac{1}{{\rm det}}\}_{\Pi}$ and $\d Y_{i,j}=0$ for all $\d\in \Lambda$ and all 
 $i,j$. We can identify $K\otimes _k R_0$ with  $K\otimes _{k_0}R_1$ where  
$R_1:=k_0\{Y_{i,j},\frac{1}{{\rm det}}\}_{\Pi}$. Let $P$ be a maximal 
$\Delta$-ideal of $R_0$.  $P$ generates an ideal  
in $K\otimes _k R_0$ which is denoted by $(P)$. Since 
$K\otimes R_0/(P)\cong K\otimes (R_0/P)\neq 0$, the ideal  
$(P)$ is a proper differential ideal. Define the ideal  
$\tilde{P}\subset R_1$ by $\tilde{P}=(P)\cap R_1$. 
By Lemma~\ref{lem2} the ideal $(P)$ is  
generated by $\tilde{P}$. If $M$ is a maximal $\Pi$-ideal of 
$R_1$ containing $\tilde{P}$ then $R_1/M$ is  simple, finitely generated $\Pi$-extension of
$C_k^\Lambda$ and so is a
constrained extension of $C_k^\Lambda$. Since $C_k^\Lambda$ is differentially closed, 
Proposition~\ref{constrainedfacts}.3 implies that 
 $R_1/M = k_0$. The  
corresponding homomorphism of $C_k^\Lambda$-algebras $R_1\rightarrow k_0$ extends to 
a differential homomorphism of $K$-algebras $K\otimes _{k_0}R_1\rightarrow K$. 
Its kernel contains $(P)\subset K\otimes _k R_0= K \otimes _{C_k^\Lambda}R_1$. 
Thus we have found a $k$-linear differential homomorphism  
$\psi :R_0\rightarrow K$ with $P\subset \ker(\psi )$. The kernel 
of $\psi$ is a differential ideal and so $P=\ker(\psi )$. 
The subring $\psi (R_0)\subset K$ is isomorphic to $R_0/P$ and 
is therefore a PPV-ring. The matrix $(\psi (X_{i,j}))$ is a fundamental  
matrix in ${\rm GL}_n(K)$ and must have the form $Z\cdot (c_{i,j})$ with 
$(c_{i,j})\in {\rm GL}_n(C_k^\Lambda)$, because the field of $\Lambda$-constants of $K$ is $C_k^\Lambda$. 
Therefore, $k\{Z, 1/\det\}_\Delta$ is a PPV-extension of $k$.    \\[0.1in]
2.~Let $K'$ be a PPV-extension of $k$ for $\d_0 Y = AY$. Part 1.~of this proposition implies that
both $K'$ and $K$ are quotient fields of PPV-rings for this equation.   Proposition~\ref{propPPVquotient} implies that there is a 
these PPV-rings are $k$-$\Delta$-isomorphic and the conclusion follows. \hspace{\fill}~$\square$\\[0.2in]
 The following result was used in Proposition~\ref{propint}.
\begin{lem} \label{PV=PPV} Let $\Delta = \{\d_0, \d_1, \ldots , \d_m \}$ and $\Lambda = \{\d_0\}$. Let \begin{eqnarray}
\d_0 Y&=& AY \label{syseqn}\\
\d_1 Y  &=&  A_1 Y\notag\\
   &\vdots&  \notag\\
  \d_m Y  &= & A_m Y \notag
  \end{eqnarray}
  be an integrable system with $A_i \in \gl_n(k)$. If $K$ is a PV-extension of $k$ for
(\ref{syseqn}), then $K$ is a PPV-extension of $k$ for $\d_0Y = AY$
 \end{lem}
\begin{proof} We first note that $C_k^\Delta$ is a subfield of $C_k^\Lambda$. Since this latter field is differentially closed, it is algebraically closed.
Therefore, $C_k^\Delta$ is also algebraically closed. The usual Picard-Vessiot theory\footnote{Proposition 1.22 of \cite{PuSi2003} proves this only for the ordinary case.  Proposition~\ref{propPPVfield} above
yields this result if we let $\Lambda = \Delta$.} implies that 
 $K$ is the quotient field of the Picard-Vessiot ring $k\{Z,1/\det Z\}_\Delta$ where
$Z$ satisfies the system (\ref{syseqn}). Since $R$ is a simple $\Delta$-ring, we have that $Z$ is
constrained over $k$, Proposition~\ref{constrainedfacts}.1 implies that  $K$ is a $\Delta$-constrained extension of $k$.  
Since $C_k^\Lambda$ is differentially closed,
Corollary~\ref{constants} implies that $C_K^{\d_0} = C_k^{\d_0}$ so $K$ is a PPV-extension of $k$.
\end{proof}

\subsection{Galois groups}\label{subgaloisgroups} In this section we shall show that the PPV-group $\Gal(K/k)$ of a
PPV-extension $K$ of $k$
is a linear differential
algebraic group and also show the correspondence between Kolchin-closed subgroups of
$\Gal(K/k)$ and $\Delta$-subfields of $K$ containing $k$.  This is done in the next Proposition
and conclusions 2.~and 3.~of Theorem~\ref{PPVthm} are immediate consequences.  \\[0.1in]
To make things a little more precise, we will use a little of the language of affine differential
algebraic geometry (see \cite{cassidy1} or \cite{kolchin_groups} for more details). We begin with some
definitions that are
the obvious differential counterparts of the usual definitions in affine algebraic geometry. Let
$k$ be a $\Delta-$field. An {\em affine differential variety} $V$ defined over $k$ is
given by a radical differential ideal $I \subset k\{Y_1, \dots , Y_n\}_\Delta$. In this case, we
shall say $V$ is a differential subvariety of affine $n$-space and write $V \subset \AX^n$.
We will identify $V$ with its {\em coordinate ring} 
$ k\{V\} = k\{Y_1, \dots , Y_n\}_\Delta/I$. Conversely, given a reduced $\Delta$-ring $R$ that is
finitely generated (in the differential sense) as a $k$-algebra, we may associate with it the
differential variety $V$ defined by the radical ideal $I$ where $R =  k\{Y_1, \dots ,
Y_n\}_\Delta/I$. 
Given any $\Delta$-field $K
\supset k$, the set of
$K$-points of $V$, denoted by $V(K)$, is the set of points of $K^n$ that are zeroes of the defining
ideal of $V$, and
may be identified with the set of $k$-$\Delta$-homomorphisms of
$k\{V\}$ to $K$.  If $V \subset \AX^n$ and $W\subset \AX^p$ are affine differential varieties defined
over $k$, a {\em differential
polynomial map} $f:V\rightarrow W$ is given by a $p$-tuple $(f_1, \ldots ,f_p) \in (k\{Y_1, \dots ,
Y_n\}_\Delta)^p$
such that the map that sends an $F \in k\{Y_1, \dots ,
Y_p\}_\Delta$ to $F(f_1, \ldots , f_p) \in k\{Y_1, \dots ,
Y_n\}_\Delta$ induces a $k$-$\Delta$-homomorphism $f^*$ of $k\{W\}$ to $k\{V\}$.  A useful
criterion for showing that a $p$-tuple $(f_1, \ldots ,f_p) \in (k\{Y_1, \dots ,
Y_n\}_\Delta)^p$ defines a differential polynomial map from $V$ to $W$ is the following: {\em 
$(f_1, \ldots ,f_p)$ defines a differential polynomial map from $V$ to $W$ if and only if for any
$\Delta$-field $K \supset k$ and any $v \in V(K)$, we have $(f_1(v), \ldots ,f_p(v)) \in W(K)$}.
This  is an easy consequence of the {\em theorem of zeros} (\cite{DAAG}, Ch.~IV.2) which in turn
is an easy consequence of the fact that a radical differential ideal is the intersection of prime
differential ideals.\\[0.1in]
Given
affine differential  varieties $V$ and $W$ defined over $k$,
we define the {\em product 
 $V\times_k W$ of $V$ and $W$} to  be the differential affine variety associated with
$k\{V\}\otimes_k
k\{W\}$. Note that since our fields have characteristic zero, this latter ring is reduced.\\
[0.1in]
In this setting, a linear differential algebraic group $G$ (defined over $k$)  is the affine
differential algebraic variety
associated with a radical differential ideal $I \subset k\{Y_{1,1,}, \ldots , Y_{n,n}, Z\}_\Delta$ such that 
\begin{enumerate}
\item $1-Z\cdot
\det((Y_{i,j})) \in I$, 
\item $({\rm id}, 1) \in G(k)$ where ${\rm id}$ is the $n\times n$ identity matrix.
\item  the map given by matrix multiplication $(g,(\det g)^{-1})(h, (\det h)^{-1}) \mapsto
(gh,(\det(gh))^{-1})$ (which is obviously
a differential polynomial map) is a map from $G\times G$ to $G$ and the inverse map $(g,\det g)^{-1})\mapsto (g^{-1},\det g)$ (also a differential polynomial map) is a map from $G$ to $G$.
\end{enumerate}
Since we assume that $1-Z\cdot
\det((Y_{i,j})) \in I$, we may assume that $G$ is defined by a radical differential ideal in the
ring $k\{Y_{1,1,}, \ldots , Y_{n,n},1/\det (Y_{i,j})\}_\Delta$, which we abbreviate as
$k\{Y, 1/\det Y\}_\Delta$. In this way, for any $K\supset k$ we may identify $G(K)$ with 
elements of $\GL_n(K)$ and the multiplication and inversion is given by the usual operations on
matrices.  We also note
that the usual Hopf algebra definition of a linear algebraic group carries over to this setting as
well. See \cite{cassidy2} for a discussion of $k$-differential Hopf algebras, and citeria for an 
affine differential algebraic group to be linear.
\begin {prop}\label{galois}                                         
Let $K\supset k$ be a PPV-field with differential Galois group 
$\PGal(K/k)$. 
Then\begin{enumerate}
\item $\PGal(K/k)$ is the group of ${C_k^\Lambda}$-points  $G(C_k^\Lambda) \subset\GL_n({C_k^\Lambda})$ of a
linear  $\Pi$-differential algebraic group $G$ over ${C_k^\Lambda}$. 
\item Let $H$ be a subgroup of $\PGal(K/k)$ satisfying $K^H = k$.  Then the Kolchin
closure $\bar H$ of $H$ is $\PGal(K/k)$.
 \item The field $K^{\PGal(K/k)}$ of
$\PGal(K/k)$-invariant elements of the Picard-Vessiot field $K$ is equal to $k$.
\end{enumerate}
\end{prop}

\begin{proof} 1.~We shall show that there is a radical $\Pi$-ideal $I \subset
S={k_0}\{Y_{i,j},\frac{1}{{\rm det}}\}_{\Pi}$ such that $S/I$ is the coordinate ring
of a linear $\Pi$-differential algebraic group $G$ and  $\Gal(K/k)$ corresponds to $G(C_k^\Lambda)$.\\[0.1in]
 Let
$K$ be the PPV-extension for the  integrable system (\ref{intsys2}). Once again we denote by $
k\{X_{i,j},\frac{1}{{\rm det}}\}_{\Pi}$ the $\Pi$-differential
polynomial ring with  the added $\Delta$-structure defined by $(\d_rX_{i,j}) = A_r(X_{i,j})$ for $\d_r \in \Lambda$.
$K$ is the field of fractions of 
$R:=k\{X_{i,j},\frac{1}{{\rm det}}\}_{\Pi}/q$, where $q$ is a maximal 
$\Delta$-ideal.
Let $r_{i,j}$ be the image of $X_{i,j}$ in $R$ so $(r_{i,j})$ is 
a fundamental matrix for the  equation $\d_0Y = AY$.  Consider the
following rings:
\[ k\{X_{i,j},\frac{1}{{\rm det}}\}_{\Pi} \subset 
K\{X_{i,j},\frac{1}{{\rm det}}\}_{\Pi} =
K\{Y_{i,j},\frac{1}{{\rm det}}\}_{\Pi} \supset 
{k_0}\{Y_{i,j},\frac{1}{{\rm det}}\}_{\Pi}\]
where the indeterminates $Y_{i,j}$ are defined by $(X_{i,j}) = (r_{i,j})(Y_{i,j})$.  Note
that $\d Y_{i,j} = 0$ for all $\d \in \Pi$. 
Since all fields are of characteristic zero, the ideal $qK\{Y_{i,j},\frac{1}{{\rm det}}\}_{\Pi}\subset
K\{X_{i,j},\frac{1}{{\rm det}}\}_{\Pi} = K\{Y_{i,j},\frac{1}{{\rm det}}\}_{\Pi}$ is a radical
$\Delta$-ideal ({\em cf.,} \cite{PuSi2003}, Corollary A.16).  Lemma~\ref{lem2} implies that
$qL[Y_{i,j},\frac{1}{{\rm det}}]$ is
 generated by $I = qK\{Y_{i,j},\frac{1}{{\rm det}}\}_{\Pi}\cap 
{k_0}\{Y_{i,j},\frac{1}{{\rm det}}\}_{\Pi}$.  Clearly $I$ is a radical $\Delta$-ideal of $S =
{k_0}\{Y_{i,j},\frac{1}{{\rm det}}\}_{\Pi}$.  We shall show that  $S/I$ is the
${\Pi}$-coordinate ring of a linear differential algebraic group, inheriting its group structure from
$\GL_n$.  In particular, we shall show that $G(C_k^\Lambda)$ is a subgroup of
$\GL_n(C_k^\Lambda)$ and that there is  an isomorphism of
$\Gal(K/k)$ onto $G(C_k^\Lambda)$.\\[0.1in]   
$\Gal(K/k)$ can be identified with the set of $(c_{i,j})
\in \GL_n(C_k^\Lambda)$ such that the  map $(X_{i,j}) \mapsto (X_{i,j})( c_{i,j})$
leaves the ideal $q$ invariant.  One can easily show that the following
statements are equivalent. 
\begin{quotation}  
\noindent (i) $(c_{i,j}) \in \Gal(K/k)$\\[0.1in] 
(ii) the map $k\{X_{i,j},\frac{1}{{\rm det}}\}_{\Pi} \rightarrow K$
defined by  $(X_{i,j}) \mapsto (r_{i,j})( c_{i,j})$ maps all elements of $q$
to zero.\\[0.1in]
 (iii) the map $K\{X_{i,j},\frac{1}{{\rm det}}\}_{\Pi}\rightarrow K$
defined by  $(X_{i,j}) \mapsto (r_{i,j})( c_{i,j})$ maps all elements of 
$qK\{X_{i,j},\frac{1}{{\rm det}}\}_{\Pi}=qK\{Y_{i,j},\frac{1}{{\rm det}}\}_{\Pi}$
to zero.\\[0.1in] 
(iv) Considering $qK\{Y_{i,j},\frac{1}{{\rm det}}\}_{\Pi}$ as an
ideal of $K\{X_{i,j},\frac{1}{{\rm det}}\}_{\Pi}$, the map 
$K\{Y_{i,j},\frac{1}{{\rm
det}}\}_{\Pi} \rightarrow K$ defined by $(Y_{i,j}) \mapsto (c_{i,j})$ sends all
elements of $qK\{Y_{i,j},\frac{1}{{\rm det}}\}_{\Pi}$ to zero. \end{quotation}
 Since the ideal 
$qK\{Y_{i,j},\frac{1}{{\rm det}}\}_{\Pi}$ is generated by
$I$, the last statement above is equivalent to $(c_{i,j})$ being a zero of the
ideal $I$, i.e., $(c_{i,j}) \in G(C^\Lambda_k)$.  Since $\Gal(K/k)$ is a group, the set
$G(C^\Lambda_k)$ is a subgroup of $\GL_n(C^\Lambda_k)$.  Therefore $G$ is a linear differential algebraic group.\\[0.2in]
2. Assuming that $\bar H \neq \Gal$, we shall derive a contradiction.  We
shall use the notation of part 1.~above.  If $\bar H \neq \Gal$, then there
exists an element $P \in {k_0}\{Y_{i,j}, \frac{1}{\det}\}_{\Pi}$ such that $P \not \in I$
and $ P(h) = 0$ for all $h \in H$.   Lemma~\ref{lem2} implies that $P
\not \in (I) = qk\{Y_{i,j},\frac{1}{{\rm det}}\}_{\Pi}$.  Let $T = 
\{ Q \in K\{X_{i,j}, \frac{1}{\det}\}_{\Pi} \ | \ 
Q \not \in (I)
\mbox{ and } Q((r_{i,j})(h_{i,j})) = 0 \mbox{ for all }h = (h_{i,j}) \in H\}$. 
Since  $K\{X_{i,j}, \frac{1}{\det}\}_{\Pi}$=$ K\{Y_{i,j}, \frac{1}{\det}\}_{\Pi}
\supset {k_0}\{Y_{i,j}, \frac{1}{\det}\}_{\Pi} $ we have that $T \neq \{0\} $.   Any
element of $K\{X_{i,j}, \frac{1}{\det}\}_{\Pi}$  may be written as $\sum_{\alpha}
f_{\alpha} Q_{\alpha}$ where  $f_{\alpha} \in K$ and $Q_{\alpha} \in k\{X_{i,j},
\frac{1}{\det}\}_{\Pi}$.  Select $  Q =  f_{\alpha_1}Q_{\alpha_1}+ \ldots +
f_{\alpha_m}Q_{\alpha_m} \in T$ with the $f_{\alpha_i}$ all nonzero and $m$
minimal.  We may assume that $f_{\alpha_1} = 1$. 
For each $h \in H$, let $ Q^h = f_{\alpha_1}^hQ_{\alpha_1}+ \ldots +
f_{\alpha_m}^hQ_{\alpha_m}$. One sees that  $ Q^h \in T$. Since $Q - Q^h$
is shorter than $Q$ and satisfies 
$(Q-Q^h)((r_{i,j})(h_{i,j})) = 0 \mbox{ for all }h = (h_{i,j}) \in H$ we must
have that $Q - Q^h \in (I)$.  If $Q-Q^h \neq 0$ then there exists an $l \in K$
such that $Q - l(Q-Q^h)$ is shorter than $Q$.  One sees that $Q - l(Q-Q^h) \in
T$ yielding a contradiction unless $Q-Q^h = 0$.  Therefore $Q = Q^h$ for all
$h \in H$ and so the $f_{\alpha_i} \in k$.  We conclude that $Q \in k\{X_{i,j},
\frac{1}{\det}\}_{\Pi}$.  Since $Q(r_{i,j}) = 0$ we have that $Q \in q$, a
contradiction.\\[0.2in]
3. Let $a =  \frac{b}{c} \in K\backslash k $ with $b,c \in R$ and let $d = 
b\otimes c - c\otimes b \in R\otimes_k R$.  Elementary properties of tensor products 
imply that $d \neq 0$ since $b$ and $c$ are linearly independent over $C_k^\Lambda$.  The ring $R\otimes_kR$
has no nilpotent elements since the characteristic of $k$ is zero ({\em cf.,} \cite{PuSi2003},
Lemma A.16).
We define a $\Delta$-ring structure on $R\otimes_k R$ by letting $\d(r_1\otimes
r_2) = \d(r_1)\otimes r_2 + r_1\otimes \d(r_2)$ for all $\d \in \Delta$. Let $J$ be
a maximal differential ideal in the differential ring
$(R\otimes_k R)[\frac{1}{d}]$.   Consider the two obvious morphisms
$\phi_i:R \rightarrow N:= (R\otimes_k R)[\frac{1}{d}]/J$.  The images of the
$\phi_i$ are generated (over $k$) by fundamental matrices of the same matrix
differential equation. Therefore both images are equal to a certain subring $S
\subset N$ and the maps $\phi_i:R \rightarrow S$ are isomorphisms.  This
induces an element $\sigma \in G $ with $\phi_1 = \phi_2\sigma$.  The image of
$d$ in $N$ is equal to $\phi_1(b)\phi_2(c) - \phi_1(c) \phi_2(b)$.  Since the
image of $d$ in $N$ is nonzero, one finds  $\phi_1(b)\phi_2(c) \neq\phi_1(c)
\phi_2(b)$.  Therefore $\phi_2((\sigma b)c) \neq \phi_2((\sigma c)b)$
and so $(\sigma b)c \neq (\sigma c )b$.  This implies $\sigma(\frac{b}{c})
\neq \frac{b}{c}$.  \end{proof} 
We have therefore completed proof of parts 2.~and 3.~ of Theorem~\ref{completegalois}.

\subsection{PPV-Rings and Torsors} \label{ppvringsandtorsors} In this section we will prove   conclusion 4.~of
Theorem~\ref{completegalois}. As in the usual Picard-Vessiot theory, this depends on identifying the
PPV-extension ring as the coordinate ring of a torsor of the PPV-group.  
\begin{defin} Let $k$ be a $\Pi$-field and $G$ a linear differential algebraic group
defined over $k$. A {\em $G$-torsor (defined over $k$)} is an  affine differential algebraic variety
$V$ defined over $k$ together with  a differential polynomial map $f:V\times_k G \rightarrow
V\times_k
V$ (denoted by  $f:(v,g) \mapsto vg$) such that 
\begin{enumerate} 
\item  for any $\Pi$-field $K\supset k, v\in V(K), g,g_1,g_2 \in G(K)$, $v1_G = v$, $v(g_1g_2) =
(vg_1)g_2$ and
\item  the
associated homomorphism $k\{V\}\otimes_k k\{V\} \rightarrow k\{V\}\otimes_k k\{G\}$ is an
isomorphism (or equivalently, for any $K \supset k$, the map $V(K)\times G(K) \rightarrow
V(K)\times V(K)$ is a bijection.
\end{enumerate}\end{defin}
We note that $V=G$ is a torsor for $G$ over $k$ with the action given by multiplication.  This
torsor is
called the {\em trivial torsor over $k$}.
We shall use the following notation. If $V$ is a differential affine variety
defined over $k$ with coordinate ring $R = k\{V\}$
and $K\supset k$ we denote by $V_K$ the differential algebraic variety (over $K$) whose coordinate ring is
$R\otimes_kK = K\{V\}$.\\[0.1in] 
We again consider the integrable system (\ref{intsys2})  over the
$\Delta$-field $k$. The PPV-ring for this equation has the
form   $R = k\{X_{i,j},\frac{1}{{\rm det}}\}_{\Pi}/q$, where $q$ is a maximal
$\Delta$-ideal. In the following, we shall think of $q$ as only a $\Pi$-differential ideal. 
We recall that $k\{X_{i,j},\frac{1}{{\rm det}}\}_{\Pi}$ is the
coordinate ring of the  linear $\Pi$-differential algebraic  group 
${\rm GL}_n$ over $k$. Let $V$ be the affine differential algebraic variety
associated with the ring $k\{X_{i,j},\frac{1}{{\rm det}}\}_{\Pi}/q$. This  is an  
irreducible and reduced $\Pi$-Kolchin-closed subset of ${\rm GL}_n$. Let $K$ 
denote the field of fractions of $k\{X_{i,j},\frac{1}{{\rm det}}\}_{\Pi}/q$.  We have
shown in the previous section that the PPV-group $\PGal(K/k) $ of
this equation may be identified with $G(C_k^\Lambda)$, that is the ${C_k^\Lambda}$-points of a
$\Pi$-linear differential algebraic group $G$ over ${C_k^\Lambda}$.  We recall how $G$ was defined.  Consider
the following rings 
\[k\{X_{i,j},\frac{1}{{\rm det}}\}_{\Pi}\subset K\{X_{i,j},\frac{1}{{\rm det}}\}_{\Pi}= 
K\{Y_{i,j},\frac{1}{{\rm det}}\}_{\Pi}\supset {k_0}\{Y_{i,j},\frac{1}{{\rm det}}\}_{\Pi} ,\] 
where the relation between the variables $X_{i,j}$ and the variables $Y_{i,j}$ 
is given by $(X_{i,j})=(r_{i,j})(Y_{i,j})$. 
The  $r_{a,b}\in K$ are the images of $X_{a,b}$ in  
$k\{X_{i,j},\frac{1}{{\rm det}}\}_{\Pi}/q\subset K$.  In Proposition~\ref{galois} we
showed that the ideal $I = qK\{X_{i,j},\frac{1}{\det}\}_{\Pi}\cap
{k_0}\{Y_{i,j},\frac{1}{\det}\}_{\Pi}$ defines  $G$.   This observation is the key to
showing the following.

\begin{prop}\label{galoistorsor} $V$ is a $G$-torsor over $k$.
\end{prop}
\begin{proof} Let $E$ be a $\Delta$-field containing $k$. 
The group $G({C_k^\Lambda}) \subset \GL_n({C_k^\Lambda})$ is precisely the set of
matrices $(c_{i,j})$ such that the map $(X_{i,j}) \mapsto (X_{i,j})(c_{i,j})$
leaves the ideal $q$ stable. In particular, for $({ c}_{i,j}) \in
G({C_k^\Lambda}), \ ({\bar z}_{i,j}) \in V(E)$ we have that $({\bar
z}_{i,j})({c}_{i,j})
\in V(E)$.  We will first show that this map defines a 
morphism from $V\times G_k \rightarrow V$.  The map is clearly defined over $k$ 
so we need only show that for any $({\bar c}_{i,j}) \in
G(E), \ ({\bar z}_{i,j}) \in V(E)$ we have that $({\bar
z}_{i,j})({\bar c}_{i,j})
\in V(E)$.  Assume that this is not true and let $({\bar c}_{i,j}) \in
G(E), \ ({\bar z}_{i,j}) \in V(K)$ be such that $({\bar z}_{i,j})({\bar c}_{i,j})
\not\in V(E)$. Let $f$ be an element of $q$ such that $f(({\bar
z}_{i,j})({\bar c}_{i,j})) \not = 0$.  Let $\{\alpha_s\}$ be a basis of
$E$ considered as a vector space over ${C_k^\Lambda}$ and let $f(({\bar
z}_{i,j})({{C}}_{i,j}))
= \sum_{\alpha_s} \alpha_sf_{\alpha_s}(({C}_{i,j}))$ where the ${C}_{i,j}$ are indeterminates and
the $f_{\alpha_s}(({C}_{i,j})) \in C_k^\Lambda\{C_{1,1}, \ldots ,C_{n,n}\}_{\Lambda}$. By assumption (and the
fact that linear independence over constants is preserved when one goes to extension fields), we
have that there is an $\alpha_s$ such that $f_{\alpha_s}((c_{i,j})) \neq 0$.  Since $C_k^\Lambda$ is a
$\Pi$-differentially closed field, there must exist $(c_{i,j}) \in G(C_k^\Lambda)$ such that 
$f_{\alpha_s}(c_{i,j}) \not = 0$.  This contradicts the fact that $f(({\bar
z}_{i,j})({ c}_{i,j})) = 0$.\\[0.1in]
Therefore the map $(V\times_k G_k)(E) \rightarrow V(E)$
defined by $(z,g) \mapsto zg$ defines a morphism $V\times_k G_k \rightarrow
V$.  At the ring level, this isomorphism corresponds to a homomorphism of
rings
\begin{eqnarray*}
k\{X_{i,j},\frac{1}{\det}\}_{\Pi}/q \rightarrow 
k\{X_{i,j},\frac{1}{\det}\}_{\Pi}/q
\otimes_{C_k^\Lambda} {C_k^\Lambda}[Y_{i,j},\frac{1}{\det}]/I\\
\hspace{1in} \simeq 
k\{X_{i,j},\frac{1}{\det}\}_{\Pi}/q
\otimes_k(k\otimes_{C_k^\Lambda} {C_k^\Lambda}\{Y_{i,j},\frac{1}{\det}\}_{\Pi}/I )
\end{eqnarray*}
where the map is induced by $(X_{i,j}) \mapsto (r_{i,j})(Y_{i,j})$.  We
have to show that the morphism
$f: V\times_k G_k\rightarrow V\times _k V$, given by $(z,g)\mapsto (zg,z)$ is an 
isomorphism of differential algebraic varieties over $k$. In terms of rings, we have to show 
that the $k$-algebra homomorphism  
$f^*:k\{V\}\otimes _k k\{V)\}\rightarrow k\{V\}\otimes _{C_k^\Lambda} k\{G\}$ is an isomorphism. 
To do this it suffices to  
find a $\Pi$-field extension $k'$ of $k$ such that $1_{k'}\otimes _k f^*$ is an  
isomorphism. For this it suffices to find $\Pi$-field extension $k'$ of $k$ such that 
$V_{k'}$ is isomorphic to $G_{k'}$ as a  $G_{k'}$-torsor over $k'$ that is, for 
some field extension $k'\supset k$, the induced morphism of varieties over $k'$, 
namely $V_{k'}\times _{k'}G_k'\rightarrow V_{k'}$, makes $V_{k'}$ into a trivial
$G$-torsor   over ${k'}$.   \\[0.1in]
Let $k' = K$, the PPV-extension of $k$ for the differential equation. We have already shown
that $I = qK\{X_{i,j},\frac{1}{\det}\}_{\Pi}\cap {k_0}\{Y_{i,j},\frac{1}{\det}\}_{\Pi}$
and this fact implies that 
\begin{eqnarray} \label{split}
K\{V\}\!=\!K\otimes _k(k\{X_{i,j},\frac{1}{{\rm det}}\}_{\Pi}/q)\! \cong\!
 K\otimes _{C_k^\Lambda}({C_k^\Lambda}\{Y_{i,j},\frac{1}{{\rm det}}\}_{\Pi}/I)\!=
 \!K\otimes_{C_k^\Lambda}C_k^\Lambda\{G\}
\!=\!K\{G\}
\end{eqnarray}
 In other words, we found an isomorphism  
$h:V_K\cong G_K$. We still have to verify that $V_K$ as  a $G$ torsor over $K$  
is, via $h$, isomorphic to the trivial torsor $G\times _{C_k^\Lambda}G_K\rightarrow
G_K$.  To do this it is enough to verify that the following diagram is
commutative and we leave this to the reader. The coordinate ring $C_k^\Lambda\{G\}$ of the group appears in
several places. To keep track of the variables, we will write $C_k^\Lambda\{G\}$ as 
${C_k^\Lambda}\{T_{i,j},\frac{1}{\det}\}_{\Pi}/\tilde{I}$ where $\tilde{I}$ is the ideal $I$ with
the variables $Y_{i,j}$ replaced by $T_{i,j}$. 

\[ \begin{array}{ccc}
K\otimes_kk\{X_{i,j},\frac{1}{\det}\}_{\Pi}/q & \stackrel{(X_{i,j})\mapsto
(X_{i,j})(T_{i,j})}{\longrightarrow} &
K\{X_{i,j},\frac{1}{\det}\}_{\Pi'}/qK\{X_{i,j},\frac{1}{\det}\}_{\Pi}\otimes_{C_k^\Lambda}C_k^\Lambda\{G\}\\&
&\\
\scriptstyle{(X_{i,j})\mapsto (r_{i,j})(Y_{i,j})}\Big\downarrow &  &
\Big\downarrow \scriptstyle{(X_{i,j})\mapsto (r_{i,j})(Y_{i,j})}\\ & & \\
K\otimes_{C_k^\Lambda}{C_k^\Lambda}\{Y_{i,j},\frac{1}{\det}\}_{\Pi}/I & \stackrel{(Y_{i,j})\mapsto
(Y_{i,j})(T_{i,j})}{\longrightarrow} &
K\{Y_{i,j},\frac{1}{\det} \}_{\Pi}/(I)_\Pi\otimes_{C_k^\Lambda}C_k^\Lambda\{G\}\end{array}   \]\end{proof}
Using this result (and its proof), we can now finish the proof of Theorem~\ref{completegalois} by proving
conclusion 4.~of this theorem. As in the usual Picard-Vessiot  theory, the proof depends on the
following group theoretic facts. Let $G$ be a linear differential algebraic group defined over
a $\Pi$-differentially closed field $C_k^\Lambda$. For any $g \in G$ the map $\rho_g:G\rightarrow
G$ given by $\rho_g(h) = hg$ is a differential polynomial isomorphism of $G$ onto $G$ and
therefore corresponds to an isomorphism $\rho^*_g:C_k^\Lambda\{G\}\rightarrow C_k^\Lambda\{G\}$.
In this way $G$ acts on the ring $k\{G\}$. Let $H$ be a
normal linear differential algebraic subgroup of $G$.  The
following facts follow from results of \cite{cassidy1} and \cite{kolchin_groups}:
\begin{enumerate} 
\item The $G$-orbit $\{\rho^*_g(f) \ | \ g \in G(C_k^\Lambda)\}$ of any $f \in C_k^\Lambda\{G\}$ spans a finite
dimensional $C_k^\Lambda$-vector space.
\item The group $G/H$ has the structure of a linear differential algebraic group (over $C_k^\Lambda$) and
its coordinate ring $C_k^\Lambda\{G/H\}$  is isomorphic to the ring of invariants $C_k^\Lambda\{G\}^H$ . 
\item The two ring $Qt(C_k^\Lambda\{G\})^H$ and $Qt(C_k^\Lambda\{G\}^H)$ are naturally $\Pi-$isomorphic, where
$Qt(..)$ denotes the total quotient ring.
\end{enumerate}
We now can prove
\begin{prop}   Let
$K$ be a PPV-extension of $k$ with Galois group $G$ and let $H$ be a normal Kolchin-closed subgroup.
Then $K^H$ is a PPV-extension of $k$.
\end{prop}
\begin{proof} Let $K$ be the quotient field of the PPV-ring $R = k\{Z, \frac{1}{\det}\}$.  As we have already
noted ({\em cf.}, (\ref{split})), we have 
\[K\otimes_k R \cong K\otimes_{C_k^\Lambda} C_k^\Lambda\{G\}\]
that is, the torsor corresponding to $R$ becomes trivial over $K$. The group $G$ acts on 
$K\otimes_{C_k^\Lambda} C_k^\Lambda\{G\}$ by acting trivially on the left factor and via $\rho^*$ on the right
factor, or trivially on the left factor and with the Galois action on the right factor.  In this way we have that $K\otimes_k R^H \cong K\otimes_{C_k^\Lambda} C_k^\Lambda\{G\}^H = K\otimes_{C_k^\Lambda}
C_k^\Lambda\{G/H\}$ and that $K\otimes_k K^H \cong K\otimes_{C_k^\Lambda}Qt(k\{G\}^H)$ by the items enumerated
above. \\[0.1in]
We now claim that $R^H$ is  finitely generated as a ${\Delta'}$-ring over $k$, hence as a $\Delta$-ring over $k$.   
Since $C_k^\Lambda\{G/H\}$
is a
finitely generated $\Pi$-$C_k^\Lambda$-algebra, we have that there exist $f_1, \ldots , f_s \in R^H$ that
generate $K\otimes_k R^H$ as a $\Pi$-$K$-algebra. We claim that $f_1, \ldots , f_s$ generate
$R^H$ as a $\Pi$-$k$-algebra.  
Let $\calM$ be a  $k$-basis of $k\{f_1, \ldots , f_s\}_{\Pi}$.
By assumption, any element of $ f \in K\otimes_k R^H = K\otimes_k k\{f_1, \ldots , f_s\}_{\Pi}$
can be written uniquely as $f = \sum_{u\in
\calM}a_u \otimes u$ where $a_u \in K$. The Galois group $G(C_k^\Lambda)$ of $K$ over $k$ also acts on
$K\otimes_k R^H$
by acting as differential automorphisms of the left factor and trivially on the right factor. Write
$1\otimes f \in 1\otimes R^H \subset K\otimes_k R^H$ as $1\otimes f = \sum_{u\in
\calM}a_u \otimes u$ where $a_u \in K$. Applying $\sigma \in G(C_k^\Lambda)$ to $f$ we have 
$1\otimes f = \sum_{u\in
\calM}\sigma(a_u) \otimes u$.  Therefore $\sigma(a_u) = a_u$ for all $\sigma \in G(C_k^\Lambda)$.  The
parameterized Galois theory implies that $a_u \in k$ for all $u$. Therefore $f \in  k\{f_1, \ldots ,
f_s\}_{\Pi}$  Ê and so $R^H = k\{f_1, \ldots , f_s\}_{\Pi}$.\\[0.1in]
Using item 1.~in the above list, we may assume that $f_1,\ldots, f_s$ form a basis of a $G/H(C_k^\Lambda)$
invariant $C_k^\Lambda$-vector space. Let $\Theta$ be the free commutative semigroup generated by
the elements of $\Lambda$. By Theorem 1, Chapter II of \cite{DAAG} (or Lemma D.11 of \cite{PuSi2003}), there exist $\theta_1=1, \ldots , \theta_s \in \Theta$
such that \[W = (\theta_i(f_j))_{1\leq i \leq s, 1\leq j \leq s}\]
is invertible. For each $\d_i \in \Lambda$, we have that $A_i = (\d_iW) W^{-1}$ is
left invariant by the action of $G/H(k_0)$.  Therefore each  $A_i \in \gl_n(k)$. Furthermore, the  $A_i$ satisfy the 
integrability conditions. We have that $K^H$ is
generated as a $\Delta$-field over $k$ by the entries of $W$. Since the constants of $K^H$ are
$C_k^\Lambda$, we have that $K^H$ is a PPV-field for the system $\d_i Y = A_i Y, \d_i \in \Lambda$. 
\end{proof}

We can now complete the proof of conclusion 4.~of Theorem~\ref{completegalois}.  If $F= K^H$
 is left
invariant by $\PGal(K/k)$ then restriction to $F$ gives a homomorphism of $\PGal(K/k)$ to
$\PGal(F/k)$. By the previous results, the kernel of this map is $H$ so $H$ is normal in
$\PGal(K/k)$.  To show surjectivity we need to show that any $\phi \in \PGal(F/k)$ extends to a
$\tilde{\phi} \in \PGal(K/k)$.  This follows from the fact the unicity of PPV-extensions.\\[0.1in]
Now assume that $H$ is normal in $\PGal(K/k)$ and that there exists an element $\tau \in \PGal(K/k)$
such that $\tau(F) \neq F$.  The Galois group of $K$ over $\tau(F)$ is $\tau H \tau^{-1}$. Since
$F \neq \tau(F)$ we have $H \neq \tau H \tau^{-1}$, a contradiction. \\[0.1in]
The last sentence of conclusion 4.~follows from the above proposition.

\subsection{Parameterized Liouvillian Extensions} In this section we will prove
Theorem~\ref{liouvillian}. One may recast this latter result in the more general setting of the last three sections but
for simplicity we will stay with the original formulation. Let $K$ and $k$ be as in the hypotheses of this theorem. Let $\KPV 
\subset K$ be the associated PV-extension as in Proposition~\ref{PPVvsPV}. \\
[0.1in]
\underline{1 $\Rightarrow$ 2:} Assume that the  Galois group $\PGal(K/k)$ contains a solvable
subgroup
of finite index. We may assume this subgroup is Kolchin closed. Since $\PGal(K/k)$ is Zariski-dense in ${\rm Gal}_{\{\d_0\}}(\KPV/k)$, we have that
this latter
group also contains a solvable subgroup of finite index. Theorem 1.43 of \cite{PuSi2003} implies
that $\KPV$ is a liouvillian extension of $k$, that is, there is a tower of $\d_0$-fields  $k = K_0
\subset K_1 \subset \ldots \subset K_r = \KPV$ such that $K_i = K_{i-1}(t_i)$ for $i = 1, \ldots r$
where either $\d_0t_i \in K_{i-1},$ or $t_i \neq 0$ and $\d_0 t_i/t_i \in K_{i-1}$ or $t_i$ is
algebraic over $K_{i-1}$. We can therefore form a tower of $\Delta$-fields $k = \tilde{K}_0
\subset \tilde{K}_1 \subset \ldots \subset \tilde{K}_r$ by inductively defining $\tilde{K}_i =
\tilde{K}_{i-1}\langle t_i \rangle_\Delta$. Since $\KPV = K_r$, we have $K = \tilde{K}_r$ and 
so $K$ is a parameterized liouvillian extension.\\
[0.1in]
\underline{3 $\Rightarrow$ 1:} Assume that $K$ is contained in a parameterized liouvillian extension
of $k$.  We wish to show that $\KPV$ is contained in a liouvillian extension of $k$.  For this we
need the following lemma.
\begin{lem}\label{lem9.4} If $L$ is a parameterized liouvillian extension of $k$ then $L = \cup_{i
\in \NX} L_i$ where $L_{i+1} = L_{i}(\{t_{i,j}\}_{j \in \NX})$ and $\{t_{i,j}\}$ is a set of
elements such that for each $j$ either $\d_0t_{i,j} \in L_i$ or $t_{i,j} \neq 0$ and $\d_0
t_{i,j}/t_{i,j} \in L_i$ or $t_{i,j}$ is algebraic over $L_i$.
\end{lem}
\begin{proof}In this proof we shall refer to a tower of fields $\{L_i\}$as above,  as a $\d_0$-tower
for $L$. By induction on the length of the tower of $\Delta$-fields defining $L$ as a 
parameterized liouvillian extension of $k$, it is enough to show the following: {\em Let $\{L_i\}$
be a $\d_0$-tower for the $\Delta$-field  $L$ and let $L\langle t\rangle_\Delta$ be an extension of
$L$ such that $\d_0t\in L,
\d_0t/t
\in L$ or $t$ is algebraic of $L$. Then there exists a $\d_0$-liouvillian tower for $L\langle
t\rangle_\Delta$.} We shall deal with three cases.\\
[0.1in] If $t$ is algebraic over $L$, then it is algebraic over some $L_{j-1}$. We then inductively
define $\tilde{L}_i = L_i \mbox{ if } i < j, \tilde{L}_{j} = L_j(t) \mbox{ and } \tilde{L}_i =
L_i(\tilde{L}_j) \mbox{ if } i>j$. The fields $\{\tilde{L}_i\}$ are then a $\d_0$-tower
for $L\langle t\rangle_\Delta$.\\
[0.1in] Now, assume that $\d_0t  = a\in L$. Let $\Theta  =\{ \d_0^{n_0}\d_1^{n_1}\cdots
\d_m^{n_m}\}$ be the commutative semigroup generated by the derivations of $\Delta$. Note that
$L\langle t
\rangle_\Delta = L(\{\theta t\}_{\theta \in \Theta})$. For any $\theta \in \Theta$ we have $\d_0
(\theta t) =
\theta (\d_0t) = \theta(a) \in L$. We define  $\tilde{L}_i = L_i(\{\theta t \ | \ (\theta a) \in
L_{i-1}\})$.  Each $\tilde{L}_i$  contains $\tilde{L}_{i-1}$ and is an extension of $\tilde{L}$ of
the correct type. Since $a \in L$,  we have that for any $\theta \in \Theta$ there exists an $i$
such that $\theta(a) \in L_{i-1}$, so $\theta t \in \tilde{L}_i$. Therefore,
$\cup_{i\in\NX}\tilde{L}_i
= L\langle t\rangle_\Delta$ so $\{\tilde{L}\}$ is a $\d_0$-tower for $L\langle t\rangle_\Delta$.\\
[0.1in]
Finally assume that $\d_0t /t = a \in L_j \subset L$. For $\theta  = \d_0^{n_0}\d_1^{n_1}\cdots
\d_m^{n_m}\in \Theta$, we define $\ord \theta = n_0+n_1+\ldots +n_m$.  For any $\theta \in
\Theta$,
 the Leibnitz rule implies that $ \theta(at) = p_\theta + a \theta t $ where 
\[
p_\theta \in \QX[\{\theta' a\}_{\ord (\theta') \leq \ord(\theta)}, \{\theta''t\}_{\ord (\theta'') <
\ord(\theta)}] \ .
\]
Note the strict inequality in the second subscript. Let $S_\theta = \{\theta' a\}_{\ord (\theta')
\leq \ord(\theta)} \cup \{\theta''t\}_{\ord (\theta'') < \ord(\theta)}$. we define a new tower
inductively:
\begin{eqnarray*}
\tilde{L}_1 & = & L_1(t)\\
\tilde{L}_i & = & \mbox{the compositum of } L_i \mbox{ and } \tilde{L}_{i-1}(\{\theta t  \ | \
S_\theta \subset \tilde{L}_{i-1}\}) 
      \end{eqnarray*}
We now show that this is a $\d_0$-tower for $L\langle t\rangle_\Delta$. We first claim that
$\tilde{L}_i$ is an $\{\d_0\}$-extension of $\tilde{L}_{i-1}$ generated by $\d_0$-integrals or
$\d_0$-exponentials of integrals or elements algebraic over $\tilde{L}_{i-1}$. 
For $i=1$, we have that $\d_0 t/t \in L_0$ and $L_1$ is generated by such elements. For $i>1$,
assume $\theta \in \Theta$ and 
$S_\theta \subset 
\tilde{L}_{i-1}$. We then have that
\[
\d_0(\frac{\theta t}{t}) = \frac{p_\theta}{t} \in \tilde{L}_{t-i}
\]
since $t, p_\theta \in \tilde{L}_{i-1}$. Therefore $\tilde{L}_{i-1}$ is generated by the correct
type of elements.\\
[0.1in]
We now show that for any $\theta \in \Theta$ there is some $j$ such that $\theta(t) \in
\tilde{L}_j$.  We proceed by induction on 
$i =\ord(\theta)$. For $i = 0$ this is true by construction.  Assume the statement is true for
$\ord(\theta') < i$. Since there are only a finite number of such 
$\theta$, there exists an $r\in \NX$ such that $\{\theta''t\}_{\ord (\theta'') < \ord(\theta)}
\subset \tilde{L}_{r}$.  Since 
$\{\theta' a\}_{\ord (\theta') \leq \ord(\theta)}$ is a finite subset of $L$, there is an $s \in
\NX$ such that $\{\theta' a\}_{\ord (\theta') \leq \ord(\theta)} \subset L_s$. Therefore for $j >\max
(r,s),
\theta t \in \tilde{L}_j$. Therefore, $\cup_{i\in\NX}\tilde{L}_i
= L\langle t\rangle_\Delta$ so $\{\tilde{L}\}$ is a $\d_0$-tower for $L\langle t\rangle_\Delta$.
\end{proof} 
Let $L$ be a parameterized liouvillian extension of $k$ containing $K$.  Lemma~\ref{lem9.4} implies
that $\KPV$ lies in a 
$\d_0$-tower. Since $\KPV$ is finitely generated, one sees that this implies that $\KPV$ lies in a
liouvillian extension of $k$. Therefore the PV-group $\Gal(\KPV/k)$ has a solvable subgroup $H$ of
finite index. Since we can identify ${\rm Gal}_{\{d_0\}}(K/k)$ with a subgroup of 
$\Gal(\KPV/k)$, we have that $\PGal(K/k) \cap H$ is a solvable subgroup of finite index in
$\PGal(K/k)$. \hfill \QED
\addcontentsline{toc}{section}{References}
\bibliographystyle{plain}
\bibliography{refs}
\end{document}